\documentclass{birkjour}
\usepackage[T1]{fontenc}

\usepackage[%
abbrev,%
backrefs,%
]{amsrefs}

\usepackage{amssymb}

\usepackage{lmodern}
\usepackage{microtype}

  \newtheorem{theorem}{Theorem}[section]

  \newtheorem{proposition}[theorem]{Proposition}
\theoremstyle{definition}
  
  \newtheorem{problem}{Problem}
\theoremstyle{remark}

\numberwithin{equation}{section}

\newcommand{\Nset}{\mathbb{N}}
\newcommand{\Rset}{\mathbb{R}}

\newcommand{\Cset}{\mathbb{C}}
\newcommand{\eps}{\varepsilon}

\newcommand{\st}{\;:\;}

\newcommand{\abs}[1]{\lvert #1 \rvert}
\newcommand{\bigabs}[1]{\bigl\lvert #1 \bigr\rvert}

\newcommand{\dualprod}[2]{\langle #1, #2 \rangle}
\newcommand{\norm}[1]{\lVert #1 \rVert}

\newcommand{\dif}{\,\mathrm{d}}

\begin{document}

%
%
%
%
%

\title[Choquard equation]
 {A guide to the Choquard equation}

\author{Vitaly Moroz}
\address{Swansea University\\ Department of Mathematics\\ Singleton Park\\
Swansea SA2~8PP\\ Wales, United Kingdom}	
\email{V.Moroz@swansea.ac.uk}

\author{Jean Van Schaftingen}
\address{Universit\'e catholique de Louvain\\ Institut de Recherche en Math\'ematique et Phy\-sique\\ Chemin du Cyclotron 2 bte L7.01.01\\ 1348 Louvain-la-Neuve \\ Belgium}
\email{Jean.VanSchaftingen@uclouvain.be}



\date{\today}
\dedicatory{In honour of Paul Rabinowitz\\
and with gratitude for his contributions\\ to the understanding of differential equations}

\begin{abstract}
We survey old and recent results dealing with the existence and properties of solutions to the Choquard type equations
\[
 -\Delta u + V(x)u = \bigl(|x|^{-(N-\alpha)} \ast \abs{u}^p\bigr)\abs{u}^{p - 2} u \qquad \text{in \(\Rset^N\)},
\]
and some of its variants and extensions.
\end{abstract}

\maketitle

\section{Introduction}

The present review paper aims to present the state of the art
around the mathematical study of the Choquard equation
\[
 -\Delta u + u = \bigl(I_\alpha \ast \abs{u}^p\bigr)\abs{u}^{p - 2} u \qquad \text{in \(\Rset^N\)},
\]
and of its variants.
Here, \(I_\alpha\) is the Riesz potential of order \(\alpha \in (0, N)\) on the Euclidean space \(\Rset^N\) of dimension \(N \ge 1\), defined for each point \(x \in \Rset^N \setminus \{0\}\) by
\begin{align*}
I_\alpha (x) &= \frac{A_\alpha}{\abs{x}^{N - \alpha}}, &
&\text{where }&
&A_\alpha = \frac{\Gamma(\tfrac{N-\alpha}{2})}
                   {\Gamma(\tfrac{\alpha}{2})\pi^{N/2}2^{\alpha} }.
\end{align*}
The nonlinearity is described by an exponent \(p \in \Rset\).

We will begin by describing briefly some physical motivations of the problem and mention some related problems (\S\ref{sectionContext}).
We will consider variants of the Choquard equation that include the introduction of an external potential \(V\), the replacement of the Riesz potential by a more general kernel and the treatment of more general nonhomogeneous nonlinearities instead of \(\abs{u}^p\).
The description will be split into to the autonomous (\S\ref{sect-Autonomous}) and nonautonomous (\S\ref{sectionNonautonomous}) cases, depending on whether the equation is invariant under translations or not.
Appendix (\S\ref{sectionAppendix}) is devoted to the relevant standard and less standard properties of the Riesz potentials.

Although we have strived to cover as extensively as possible the existing literature on the topic, the choice of the topics that we are covering might not appear to the reader as consistent as we have intended it and we have very likely overlooked relevant works in the abundance of independent threads of literature on the Choquard equation.

\section{Context}
\label{sectionContext}

\subsection{Physical models}

The Choquard equation
\begin{equation}\label{eqChoquard3d}
 -\Delta u + u = \bigl(I_2 \ast \abs{u}^2\bigr)u \qquad \text{in \(\Rset^3\)},
\end{equation}
has
appeared in the context of various physical models.
It seems to originate from H.\thinspace{}Fr\"ohlich and S.\thinspace{}Pekar's model of the \emph{polaron}, where free electrons in an ionic lattice interact
with phonons associated to deformations of the lattice or with the polarisation that it creates on the medium (interaction of an electron with its own hole) \citelist{\cite{Pekar1954}\cite{Frohlich1937}\cite{Frohlich1954}}.
The Choquard equation was also introduced by Ph.\thinspace{}Choquard in 1976 in the modelling of a \emph{one-component plasma} \cite{Lieb1977}.

In general, the associated Schr\"odinger-type evolution equation
\begin{equation}
\label{eqHartree}
 i \partial_t \psi = \Delta \psi + (W \ast \abs{\psi}^2) \psi
\end{equation}
is a model large system of non-relativistic bosonic atoms and molecules
with an attractive interaction that is weaker and has a longer range than that of the nonlinear Schr\"odinger equation (where the interaction potential \(W\) is formally Dirac's delta at the origin) \cite{FrohlichLenzmann2004}.
The equation \eqref{eqHartree} arises as a mean-field limit of a bosonic system with attractive two-body interactions; this limit can be taken rigorously in many cases  \citelist{\cite{FrohlichLenzmann2004}\cite{LewinNamRougerie2014}}.

In showing that his polaron model arises as an asymptotic limit of the Fr\"ohlich polaron,
S.\thinspace{}Pekar had conjectured that the groundstate level of the Pekar polaron problem should be characterised in terms of \emph{Brownian motion}
\begin{multline}
 \lim_{\alpha \to \infty} \lim_{t \to \infty} \frac{1}{t} \log \biggl(\mathbb{E} \exp \Bigl(\alpha \int_0^t \int_0^t \frac{e^{-\abs{\sigma - s}}}{\abs{x (\sigma) - x (s)}} \dif \sigma \dif s \Bigr)\biggr)\\
 = \sup\, \biggl\{ 2 \int_{\Rset^3}\int_{\Rset^3} \frac{\abs{u (x)}^2\abs{u (y)}^2}{\abs{x - y}}\dif x \dif y - \int_{\Rset^3} \abs{\nabla u}^2 \st \int_{\Rset^3} \abs{u}^2 = 1 \biggr\}
\end{multline}
where \(\mathbb{E}\) is the expectation with respect to the three-dimensional Brownian motion \(x (\cdot)\) tied at both ends (paths \(x(0) = 0 = x (t)\)).
This conjecture was proved by M.\thinspace{}D.\thinspace Donsker and S.\thinspace R.\thinspace S.\thinspace Varadhan \citelist{\cite{DonskerVaradhan1981}\cite{DonskerVaradhan1983}}.
Another mathematical analysis of the asymptotics of the Fr\"ohlich polaron was provided by E.\thinspace{}H.\thinspace{}Lieb and L.\thinspace{}E.\thinspace{}Thomas \cite{LiebThomas1997}.

Finally, the Choquard equation is also known as the \emph{Schr\"o\-ding\-er--Newton equation} in models coupling the Schr\"o\-ding\-er equation of \emph{quantum physics} together with nonrelativistic \emph{Newtonian gravity}
\citelist{%
\cite{Diosi1984}%
\cite{Penrose1996}%
\cite{Penrose1998}%
\cite{Jones1995gravitational}%
\cite{Jones1995newtonian}%
\cite{Bahramietal2014}
} (see also \citelist{\cite{RuffiniBonazzola1969}\cite{Moller1962}\cite{Rosenfeld1963}} for relativistic versions).
The equation can also be derived from the Einstein--Klein--Gordon and Ein\-stein--Dirac system \cite{GiuliniGrossardt2012}.
Such a model has been proposed for boson stars \cite{SchunckMielke2003} and for the collapse of galaxy fluctuations of scalar field dark matter \citelist{\cite{GuzmanUrenaLopes2003}\cite{GuzmanUrenaLopes2004}}.
Further models have been developed including a gravitomagnetic potential \cite{Manfredi2015}
and self-field coupling \cite{Franklinetat2015}.

\subsection{Related equations}

The Choquard equation is related to several other partial differential equations with nonlocal interactions,
which will be outside of the scope of the present survey.

Standing wave solutions to focusing Hartree equation
\begin{equation}
\label{eqHartreeRiesz}
  i \partial_t \psi = \Delta \psi + \bigl(I_\alpha \ast \abs{\psi}^2\bigr) \psi
\end{equation}
are solutions to the Choquard equation.
The local existence of solutions is known (see for example \cite{GinibreVelo1980}).
The global existence is more delicate, due to the focusing character of the nonlinearity. Some results are available when the nonlinearity mildly delocalised (\(\alpha > N - 2\)) \cite{GinibreVelo1980}*{theorem~3.1}.
The global existence of solutions to
\begin{equation}
\label{eqHartreepRiesz}
  i \partial_t \psi = \Delta \psi + \bigl(I_\alpha \ast \abs{\psi}^p\bigr) \abs{\psi}^{p - 2} \psi
\end{equation}
has been studied for \(2 \le p < 2 + \frac{4}{N - 2}\) and \(\alpha = N - 2\) \cite{GenevVenkov2012}.
Blow-up and soliton dynamics have also been studied \citelist{\cite{GenevVenkov2012}\cite{BonannoDAveniaGhimentiSquassina2014}\cite{DAveniaSquassina2014}}.

Back to the stationary setting, the interaction potential instead of being attractive as in the Choquard equation can be taken to be repulsive.
In the presence of an external potential \(V \in C (\Rset^N; \Rset)\), the stationary Hartree equation with a Coulombic potential
\[
 -\Delta u + V u + \bigl(I_2 \ast \abs{u}^2\bigr)u = \lambda u
\]
has been studied by many authors (see for example \citelist{\cite{GustafsonSather1971}\cite{BenguriaBrezisLieb81}\cite{Lions1987}\cite{LeBrisLions2005}}).
Instead of (or in addition to) imposing an external potential, a local nonlinear interaction can be added.
This leads to the Schr\"odinger--Poisson (Schr\"odinger--Poisson--Slater) equation of the form
\[
 -\Delta u + V u + \bigl(I_2 \ast \abs{u}^2\bigr)u = \abs{u}^{q-2} u,
\]
which has been the object of multiple works (see for example \citelist{\cite{BenciFortunato1998}\cite{ArriolaSoler2001}\cite{Ruiz2006}\cite{AmbrosettiRuiz2008}\cite{ZhaoZhao2008}\cite{Ruiz2010}\cite{IanniRuiz2012}})
and survey papers \citelist{\cite{Ambrosetti2008}\cite{CattoDolbeaultSanches13}}.

\section{Autonomous Choquard equation}\label{sect-Autonomous}

In this section we explore the solutions to the Choquard equation in the case
where the problem is invariant under translations of the Euclidean space \(\Rset^N\).
An interesting family of the problems which extends \eqref{eqChoquard3d} is given by the autonomous homogeneous Choquard equations:
\begin{equation}
 \label{eqAutoHomoChoquard}
 -\Delta u + u = \bigl(I_\alpha \ast \abs{u}^p\bigr) \abs{u}^{p - 2} u \qquad \text{in \(\Rset^N\)},
\end{equation}
where $N\in\mathbb N$, $\alpha\in(0,N)$ and $p>1$.

\subsection{Groundstates}

Solutions of problem~\eqref{eqAutoHomoChoquard} are, at least formally,
critical points of the action functional \(\mathcal{A}\) defined for a function \(u : \Rset^N \to \Rset\) by
\begin{equation}
\label{eqAutoHomoAction}
 \mathcal{A} (u)
 = \frac{1}{2} \int_{\Rset^N} \big(\abs{\nabla u}^2 + \abs{u}^2\big)
 - \frac{1}{2 p} \int_{\Rset^N} \bigl(I_\alpha \ast \abs{u}^{p}\bigr) \abs{u}^p.
\end{equation}

\subsubsection{Functional setting}

The first term in the definition of the action functional \(\mathcal{A}\) by \eqref{eqAutoHomoAction}
suggests taking naturally as the domain the classical Sobolev space \(H^1 (\Rset^N)\) of functions in \(L^2 (\Rset^N)\) whose weak derivative is also square-integrable.
This raises then the question whether the second term,
which involves the convolution, is well-defined and sufficiently smooth.

If we assume that \(u \in L^{\frac{2N p}{N + \alpha}} (\Rset^N)\) and if we apply the Hardy--Littlewood--Sobolev inequality \eqref{eqHLS} to the function \(f = \abs{u}^p \in L^\frac{2N}{N + \alpha} (\Rset^N)\), we obtain, in view of H\"older's inequality
\begin{equation}
\begin{split}
  \int_{\Rset^N} \bigl(I_\alpha \ast \abs{u}^p\bigr) \abs{u}^p
  &\le \Bigl(\int_{\Rset^N} \bigabs{I_\alpha \ast \abs{u}^p}^\frac{2N}{N - \alpha} \Bigr)^{\frac{1}{2} - \frac{\alpha}{2 N}}
  \Bigl(\int_{\Rset^N} \abs{u}^\frac{2N p}{N + \alpha}\Bigr)^{\frac{1}{2} + \frac{\alpha}{2 N}}\\
  &\le C_{N, \alpha, 2N/(N + \alpha)}
  \Bigl(\int_{\Rset^N} \abs{u}^\frac{2 N p}{N + \alpha} \Bigr)^{1 + \frac{\alpha}{N}}.
\end{split}
\end{equation}
It remains then to determine when the condition \(u \in L^{\frac{2N p}{N + \alpha}} (\Rset^N)\)
is satisfied. By the classical Sobolev embedding theorem, there is a continuous
embedding \(H^1 (\Rset^N) \hookrightarrow L^\frac{2 N p}{N + \alpha} (\Rset^N)\) with \(r \in [1, \infty)\) if and only if
\(\frac{1}{2} - \frac{1}{N} \le \frac{N + \alpha}{2 N p} \le \frac{1}{2}\), or, equivalently,
\(\frac{N - 2}{N + \alpha} \le \frac{1}{p} \le \frac{N}{N + \alpha}\).
Moreover, we have then
\begin{equation}
 \int_{\Rset^N} \bigl(I_\alpha \ast \abs{u}^p\bigr) \abs{u}^p
 \le C \Bigr(\int_{\Rset^N} \abs{\nabla u}^2 + \abs{u}^2 \Bigr)^p,
\end{equation}
for some constant \(C\) depending on \(N\), \(\alpha\) and \(p\).

Basic differentiability properties of the functional \(\mathcal{A}\) follow
from these estimate by classical nonlinear functional analysis arguments:

\begin{proposition}
\label{propositionAutoHomogFunctional}
If \(p \in (1, \infty)\) satisfies
\[
  \frac{N - 2}{N + \alpha} \le \frac{1}{p} \le \frac{N}{N + \alpha},
\]
then the functional \(\mathcal{A}\) is well-defined and continuously
Fr\'echet--differentiable on the Sobolev space \(H^1 (\Rset^N)\).\\
If moreover \(p \ge 2\), then the functional \(\mathcal{A}\) is twice continuously
Fr\'echet--diffe\-ren\-tiable.
\end{proposition}

Proposition~\ref{propositionAutoHomogFunctional} follows from the corresponding differentiability properties of the superposition map
\(u \in H^1 (\Rset^N) \mapsto \abs{u}^p \in L^\frac{2N}{N + \alpha} (\Rset^N)\) and from the smoothness of the quadratic form involving the Riesz potential, resulting from the boundedness of the latter.

Restrictions on the exponents in nonlinearities related to Sobolev embedding theorems are a very classical feature of semilinear elliptic problems.
The appearance of
a lower nonlinear restriction in the class of Choquard-type problems is more remarkable: the lower critical exponent \(\frac{N + \alpha}{N}\)
is strictly greater than \(1\); the condition for well-definiteness is thus more stringent
that imposing a superlinearity condition \(p > 1\).

Another interesting point is that although the function is well-defined
and the nonlinearity is superlinear, twice differentiability only occurs when \(p \ge 2\).
It will appear in the sequel that some properties of the solutions can also change
dramatically when the exponent \(p\) crosses the value \(2\).

\subsubsection{Existence of solutions}
We define a solution \(u \in H^1 (\Rset^N)\) to be a \emph{groundstate}
of the Choquard equation \eqref{eqAutoHomoChoquard}
whenever it is a solution that minimises the action functional \(\mathcal{A}\)
among all nontrivial solutions.
Groundstates exists when the nonlinearity exponent \(p\) satisfies strictly the
inequalities for the well-definiteness of proposition~\ref{propositionAutoHomogFunctional}
\citelist{\cite{Lieb1977}\cite{Menzala1980}\cite{MorozVanSchaftingen2013JFA}\cite{MaZhao2010}
\cite{Lions1984-1}*{theorem~III.1}\cite{GenevVenkov2012}*{lemma 2.7}\cite{Efinger1984}*{theorem~3.1}}.

\begin{theorem}
\label{theoremExistenceGroundStates}
If
\[
 \frac{N - 2}{N + \alpha} < \frac{1}{p} < \frac{N}{N + \alpha},
\]
then \eqref{eqAutoHomoChoquard} has a groundstate solution \(u \in H^1 (\Rset^N)\).
\end{theorem}

In general, the groundstate \(u\) can be constructed by showing that the infimum
on the \emph{Nehari manifold}
\[
 \inf \,\bigl\{\mathcal{A} (u) \st u \in H^1 (\Rset^N) \setminus \{0\}
   \text{ and } \dualprod{\mathcal{A}' (u)}{u} \bigr\}
\]
is achieved. This is equivalent to prove that the mountain-pass minimax level
\citelist{\cite{AmbrosettiRabinowitz1973}\cite{Rabinowitz1986}\cite{Rabinowitz1995}}
\[
  \inf_{\gamma \in \Gamma} \sup_{[0, 1]} \mathcal{A} \circ \gamma,
\]
where the class of paths \(\Gamma\) is defined by
\[
\Gamma = \bigl\{ \gamma \in C([0, 1]; H^1 (\Rset^N))\st
  \gamma (0) = 0 \text{ and } \mathcal{A} (\gamma (1)) < 0\bigr\},
\]
is a critical value.
The minimisers of the Sobolev-like quotient
\begin{equation}\label{e-min-Sobolev}
  \frac{\displaystyle\int_{\Rset^N} \abs{\nabla u}^2 + \abs{u}^2}
    {\Bigl(\displaystyle \int_{\Rset^N} \bigl(I_\alpha \ast \abs{u}^p\bigr) \abs{u}^p\Bigr)^\frac{1}{p}}
\end{equation}
are groundstates up to multiplication by a constant
and the minimisers of the Weinstein-like \cite{Weinstein1982} quotient
\[
 \frac{\Bigl(\displaystyle\int_{\Rset^N} \abs{\nabla u}^2\Bigr)^{\frac{N}{2} - \frac{N + \alpha}{2 p}}
 \Bigl(\displaystyle\int_{\Rset^N} \abs{u}^2\Bigr)^{\frac{N + \alpha}{2 p} - \frac{N - 2}{2}}}
 {\Bigl(\displaystyle\int_{\Rset^N} \bigl(I_\alpha \ast \abs{u}^p\bigr) \abs{u}^p\Bigr)^\frac{1}{p}},
\]
are groundstates of \eqref{eqAutoHomoChoquard}  to multiplication and rescaling in space
\citelist{\cite{MorozVanSchaftingen2013JFA}*{proposition~2.1}\cite{GenevVenkov2012}}.

Under the additional assumption that \(p < 1 + \frac{\alpha + 2}{N}\),
minimisers of
\[
 \inf \Bigl\{\frac{1}{2} \int_{\Rset^N} \abs{\nabla u}^2
   - \frac{1}{2p} \int_{\Rset^N} \bigl(I_\alpha \ast \abs{u}^p\bigr)\abs{u}^p
 \st u \in H^1 (\Rset^N) \text{ and } \int_{\Rset^N} \abs{u}^2 \le \lambda
 \Bigr\}
\]
are groundstates of \eqref{eqAutoHomoChoquard} up to multiplication and rescaling in the Euclidean space \(\Rset^N\)
\cite{MorozVanSchaftingen2013JFA}*{proposition~2.1}; this latter
formulation was the original framework
\citelist{\cite{Lieb1977}\cite{Lions1984-1}\cite{Efinger1984}\cite{FrohlichLenzmann2004}}.

Since the equation \eqref{eqAutoHomoChoquard}
and the associated functional \(\mathcal{A}\) are invariant under translations,
the variational formulation lacks compactness properties that would make it straightforward.

A first way to handle this compactness issue is to rely on inequalities for
symmetrisation by rearrangement
(see for example \citelist{\cite{PolyaSzego1951}\cite{Kawohl1985}
\cite{BrockSolynin2000}\cite{LiebLoss2001}*{chapter 3}})
to restrict without loss of generality
the minimisation argument to radial functions
\citelist{\cite{Lieb1977}\cite{Menzala1980}\cite{GenevVenkov2012}*{lemma 2.7}\cite{Efinger1984}}.

A more delicate but more robust approach is to apply the concentration compactness
method of P.-L. Lions \citelist{\cite{Lions1984-1}\cite{FrohlichLenzmann2004}\cite{MorozVanSchaftingen2013JFA}}.
It is useful to note that a counterpart of the
Brezis--Lieb \cite{BrezisLieb1983} lemma holds for Riesz potentials
\citelist{%
\cite{Ackermann2006}*{\S 5.1}%
\cite{MorozVanSchaftingen2013JFA}%
\cite{BellazziniFrankVisciglia}%
\cite{YangWei2013}*{lemma 3.2}%
\cite{MercuriMorozVanSchaftingen}%
}: if the sequence \((u_n)_{n \in \Nset}\) converges weakly to \(u\) in \(H^1 (\Rset^N)\), then
\begin{multline}
\label{eqBrezisLieb}
 \lim_{n \to \infty} \int_{\Rset^N} \bigl(I_\alpha \ast \abs{u_n}^p\bigr) \abs{u_n}^p -
 \int_{\Rset^N} \bigl(I_\alpha \ast \abs{u - u_n}^p\bigr) \abs{u - u_n}^p\\
 = \int_{\Rset^N} \bigl(I_\alpha \ast \abs{u}^p\bigr) \abs{u}^p.
\end{multline}

\subsection{Other solutions}


Radial solutions of the Choquard equation \eqref{eqAutoHomoChoquard} have also been constructed when \(p = 2\) and \(\alpha = 2\)
by \emph{variational methods} \cite{Stuart1980}*{\S 5}
and by \emph{shooting methods} for systems of ordinary differential equations \citelist{\cite{ChoquardStubbeVuffray2008}%
\cite{TodMoroz1999}%
}.

Besides the groundstate, variational methods can be used to show
the existence of many other solutions.

The existence of infinitely many solutions has been proved when \(p = 2\) \citelist{\cite{Bongers1980}\cite{Lions1980}} (see also \cite{Lions1982}).
The construction relies on the \emph{Krasnosel\cprime{}ski\u\i{} genus} \cite{Krasnoselskii1964}
and on \emph{associated minimax theorems of Ambrosetti and Rabinowitz}
\cite{AmbrosettiRabinowitz1973}*{theorem~2.8} (see also \citelist{\cite{Rabinowitz1974}\cite{Rabinowitz1986}\cite{Rabinowitz1995}\cite{MawhinWillem1989}*{theorem~6.1}\cite{Willem1996}*{chapter 3}\cite{Struwe2008}*{theorem~5.7}}) applied to the subspace of radial functions.
Moreover, these radial solutions can have an arbitrary number of nodal domains \cite{Weth2001}*{theorem~9.5}.
More generally, there exist infinitely many solutions satisfying \(u \circ g = \tau (g) u\) for every \(g\) in a group of linear isometries \(G \subseteq O (N)\) such that every \(x \in \Rset^N \setminus \{0\}\) has an infinite orbit and a homomorphism \(\tau : G \to \{-1, 1\}\)  \cite{CingolaniClappSecchi2012}. Under some weaker condition on the orbits, there are one or several invariant solutions \cite{ClapSalazar2013}.
Most of these results and their proof are valid in the whole intercritical range \(\frac{N - 2}{N + \alpha} < \frac{1}{p} < \frac{N}{N + \alpha}\), although the proof are written for narrower ranges.

When \(\frac{N - 2}{N + \alpha} < \frac{1}{p} < \frac{N}{N + \alpha}\) there are also solutions that \emph{minimise the action functional among solutions that are odd} with respect to a hyperplane \cite{GhimentiVanSchaftingen}.
When moreover  \(p \ge 2\), there exists a solution that \emph{minimises the action functional \(\mathcal{A}\)
among all sign-changing solutions} \citelist{\cite{GhimentiVanSchaftingen}\cite{GhimentiMorozVanSchaftingen}}.
When \(\alpha\) is either close to \(0\) or close to \(N\), then the minimal nodal solution is odd \cite{RuizVanSchaftingen2016}.

\begin{problem}
Is the minimal nodal solution odd for \emph{each} $\alpha\in(0,N)$?
\end{problem}

The solutions are constructed by minimising the action on the \emph{Nehari nodal set}
\citelist{\cite{CeramiSoliminiStruwe1986}\cite{CastroCossioNeuberger1997}\cite{CastorCossioNeuberger1998}}. Remarkably, such solutions cannot appear in the autonomous case for the local nonlinear Schr\"odinger equation.
A proof of the existence of minimal action radial nodal solutions has been proposed \cite{Ye2015}; unfortunately that proof does not seem to address the question whether the constructed solution changes sign.

\subsection{Properties of solutions}
Solutions of the Choquard equation \eqref{eqAutoHomoChoquard} enjoy various qualitative properties: weak solutions turn out to be classical solutions, groundstates are, up to translation and inversion of the sign, positive and radially symmetric functions, the asymptotic rate of decay at infinity is described precisely, they are sometimes known to be unique and nondegenerate, they satisfy the Poho\v{z}aev variational identity.

\subsubsection{Regularity of solutions}\label{sectRegularity}
The solutions constructed variationally are weak solutions to the Choquard equation \eqref{eqAutoHomoChoquard}: for every test function \(\varphi \in H^1 (\Rset^N)\),
\[
 \int_{\Rset^N} \nabla u \cdot \nabla \varphi + u \varphi = \int_{\Rset^N} \bigl(I_\alpha \ast \abs{u}^p\bigr) \abs{u}^{p - 2}u \varphi.
\]
The classical bootstrap method for subcritical semilinear elliptic problems combined
with estimates for Riesz potentials allows to prove that \emph{any weak solution is a function of class \(C^2\)} (twice continuously differentiable); if moreover the solution \(u\) is positive or \(p\) is an even integer, then \(u\) is of class \(C^\infty\) \citelist{\cite{Lieb1977}*{theorem~8}
\cite{Menzala1980}*{theorem~6.2}\cite{CingolaniClappSecchi2013}*{lemma A.1}\cite{MorozVanSchaftingen2013JFA}*{proposition~4.1}
\cite{Lei2013MZ}\cite{Lei2013SIAN}} (see also works on a related nonlinear integral equation \cite{ChenLiOu2005}).
A nonlocal counterpart of the regularity result of Brezis and Kato for local elliptic operators \cite{BrezisKato1979}*{theorem~2.3} (see also \cite{Trudinger1968}*{theorem~3}) can be used to cover the \emph{critical exponent} \(p = \frac{N + \alpha}{N - 2}\) \cite{MorozVanSchaftingen2015TAMS}*{theorem~2}.

The behaviour of positive subsolutions \(u \in C^2 (\Rset^N \setminus \{0\}) \cap L^p (\Rset^N)\) of the equation
\[
 -\Delta u = \bigl(I_\alpha \ast u^p\bigr) u^q \qquad \text{in \(B_1 \setminus \{0\}\)}
\]
has been studied.
Depending on the values of the parameters, the singularity at \(0\) can be  removable, can behave like \(1/\abs{x}^{N - 2}\) or have arbitrary growth \cite{GherguTaliafero} (see also \cite{ChenZou} for related results).
In the supercritical case \(\alpha = 2\) and \(p = \frac{2N}{N - 2}\), a blow-up analysis of solutions of \eqref{eqAutoHomoChoquard} has been performed in terms of the behaviour of the Riesz-potential term \(I_\alpha \ast \abs{u}^p\) \cite{Disconzi2013}, however the existence of solutions in this supercritical regime is not yet known.

\subsubsection{Positivity of groundstates}

For every \(u \in H^1 (\Rset^N)\), one has \(\mathcal{A} (\abs{u}) = \mathcal{A} (u)\).
Therefore if \(u\) is a groundstate of the Choquard equation \eqref{eqAutoHomoChoquard}, then the function \(\abs{u}\) is also a solution. By the regularity properties of solutions (\S\ref{sectRegularity}) and by the strong maximum principle for second-order differential operators, the function \(u\) is continuous and \(\abs{u} > 0\); therefore any groundstate \(u\) cannot vanish anywhere in \(\Rset^N\) and has thus \emph{constant sign} on the whole \(\Rset^N\) \cite{MorozVanSchaftingen2013JFA}*{proposition~5.1}.
If the equation \eqref{eqAutoHomoChoquard} is considered with \emph{complex valued functions} \(u\), then \(u\) does not vanish and has constant phase \cite{CingolaniSecchiSquassina2010}*{lemma 2.10}.

\subsubsection{Symmetry of groundstates and positive solutions}\label{sectSymmetry}

In the constructions of a groundstate relying on a \emph{symmetrisation}6 argument,
the groundstate is by design radially symmetric and nonincreasing \citelist{\cite{Lieb1977}\cite{Menzala1980}\cite{GenevVenkov2012}*{lemma 2.7}\cite{Efinger1984}}.
If \(u^*\) is the symmetric decreasing rearrangement of the function \(u \in H^1 (\Rset^N)\), in view of the cases of equality in the Riesz--Sobolev convolution inequality \eqref{eqRearrangement},  \(\mathcal{A} (u) = \mathcal{A} (u^*)\) if and only if \(u\) is the translation of radially symmetric and nonincreasing function \citelist{\cite{Lieb1977}*{lemma 3}\cite{Burchard1996}}, and therefore
\emph{any groundstate is symmetric} \cite{Lieb1977}.
The approach works in fact better for Choquard problems than for local problems, where only the
much less conclusive equality cases of the P\'olya--Szeg\H o inequality \cite{BrothersZiemer1988}
are available.

\emph{Polarisations,} also called
two-point rearrangements \citelist{\cite{Baernstein1994}\cite{BrockSolynin2000}}, have been used to prove the symmetry of groundstates for some local problems \citelist{\cite{BartschWethWillem2005}\cite{VanSchaftingenWillem2008}}.
The approach turns out to be even simpler for Choquard problems \cite{MorozVanSchaftingen2013JFA} (see also \cite{MorozVanSchaftingen2015TAMS}).

When \(\alpha = 2\), the symmetry of groundstates has also been proved by showing that one of the extensions by even reflection of the restrictions of the groundstate to two halfspaces separated by a given hyperplane is also a groundstate, implying symmetry by a suitable \emph{unique continuation principle} \cite{LopesMaris2008}.

The symmetry of \emph{positive solutions} has also been studied by the \emph{moving plane
method} \cite{GidasNiNirenberg1979}.
For the Choquard equation, the symmetry of positive solutions has been proved when the parameters of the equation and some intermediate exponents satisfy a set of inequalities \cite{MaZhao2010}; a sufficient condition given for the symmetry of positive solutions is that \(p\ge 2\) and
\begin{equation}
\label{conditionMaZhao}
 [2, \tfrac{2 N}{N - 2}]
  \cap (p, \tfrac{p N}{\alpha})
  \cap (\tfrac{(2 p - 2) N}{\alpha + 2}, \tfrac{(2 p - 1) N}{\alpha + 2})
  \cap [\tfrac{(2 p - 1) N}{N + \alpha}, \infty)
 \ne \emptyset.
\end{equation}
The question whether these conditions are always satisfied in the case where symmetrisation and polarisations argument work has yet to be clarified.
The moving plane method has been applied also to other problems with a nonlocal nonlinearity
\citelist{\cite{FrankLenzmann}\cite{ChenLiOu2006}}.

\subsubsection{Decay of groundstates}%
\label{s-Decay}%

The behaviour of a groundstate \(u\) to the homogeneous Choquard equation \eqref{eqAutoHomoChoquard}
has been studied in detail \cite{MorozVanSchaftingen2013JFA}:
when either \(p > 2\) or \(p = 2\) and \(\alpha > N - 1\),
then the decay is exponential and there exists \(c \in (0, \infty)\), depending on \(u\), such that
\begin{equation}
 u (x) = \bigl(c + o (1)\bigr)\frac{e^{-\abs{x}}}{\abs{x}^\frac{N - 1}{2}} \qquad \text{as \(\abs{x} \to \infty\)};
\end{equation}
when \(p = 2\) and \(\alpha \in [N - 1, N)\), then the decay is a mild perturbation of the previous one: there exists \(c \in (0, \infty)\) such that
\begin{equation}
 u (x) =  \bigl(c + o (1)\bigr)\frac{e^{-\abs{x}}}{\abs{x}^{\frac{N - 1}{2}}} \exp \int_{\nu}^{\abs{x}} \Bigl(1 - \sqrt{1 - \tfrac{\nu^{N - \alpha}}{s^{N - \alpha}}}\,\Bigr) \dif s  \qquad \text{as \(\abs{x} \to \infty\)},
\end{equation}
where
\begin{equation}
  \nu^{N - \alpha} = \frac{\Gamma(\tfrac{N-\alpha}{2})}{\Gamma(\tfrac{\alpha}{2})\pi^{N/2}2^{\alpha}} \int_{\Rset^N} \abs{u}^2
  = \frac{\Gamma(\tfrac{N-\alpha}{2})}{\Gamma(\tfrac{\alpha}{2})\pi^{N/2}2^{\alpha}}(\alpha + 4 - N)\mathcal A(u);
\end{equation}
when finally \(p \in (1 + \frac{\alpha}{N}, 2)\) the decay is polynomial and
\begin{equation}
  u (x) = \bigl(1 + o (1)\bigr) \frac{1}{\abs{x}^\frac{N - \alpha}{2 - p}} \Bigl(\frac{\Gamma(\tfrac{N-\alpha}{2})}{\Gamma(\tfrac{\alpha}{2})\pi^{N/2}2^{\alpha}} \int_{\Rset^N} \abs{u}^p\Bigr)^\frac{1}{2 - p}
 \qquad \text{as \(\abs{x} \to \infty\)}.
\end{equation}
Quite remarkably, the groundstates are always \emph{localised} enough (in the sense of \eqref{eqRieszLocalised}) so that
\begin{equation}
 \bigl(I_\alpha \ast \abs{u}^p\bigr) (x) = \Bigl(I_\alpha(x) \int_{\Rset^N} \abs{u}^p\Bigr) \bigl(1 + o (1)\bigr), \qquad \text{as \(\abs{x} \to \infty\)};
\end{equation}
these asymptotic bounds coincide with the asymptotic bounds for supersolutions \cite{MorozVanSchaftingen2013JDE}.

\subsubsection{Uniqueness and nondegeneracy}\label{s-nodegeneracy}

When \(\alpha = 2\) and \(p = 2\), there exists \emph{at most one} radial positive solution to the Choquard equations \citelist{\cite{Lieb1977}\cite{TodMoroz1999}\cite{Lenzmann2009}\cite{MaZhao2010}\cite{WangYia2016}}.

In general, the set of groundstates is known to be compact up to translations
\citelist{%
\cite{CingolaniSecchiSquassina2010}*{proposition~2.14}%
\cite{MorozVanSchaftingen2015TAMS}*{corollary 4.2}%
}.

When \(N = 3\), \(\alpha = 2\) and \(p = 2\), this groundstate is \emph{nondegenerate} up to translations \citelist{\cite{Lenzmann2009}\cite{WeiWinter2009}}, that is, the kernel of the linearised operator \(\mathcal{A}'' (u)\) is generated by the directional derivatives of the solution \(u\): if \(\varphi \in H^1 (\Rset^N)\)
\begin{equation}
\label{eqNondegeneracy}
 -\Delta \varphi + \varphi = (p-1) \bigl(I_\alpha \ast \abs{u}^p\bigr) \abs{u}^{p - 2} \varphi
 + p \bigl(I_\alpha \ast (\abs{u}^{p - 2}u \varphi) \bigr) \abs{u}^{p - 2} u,
\end{equation}
then there exists \(\xi \in \Rset^N\) such that \(\varphi = \xi \cdot \nabla u\).

This nondegeneracy allows to extend the uniqueness and nondegeneracy to the slightly superquadratic case \(N = 3\), \(\alpha = 2\) and \(p > 2\) close to \(2\) \cite{Xiang}.

\subsubsection{Poho{\v{z}}aev identity}
If \(N \in \Nset\), \(\alpha \in (0, N)\) and \(p > 1\), any solution \(u \in W^{1, 2} (\Rset^N) \cap L^{\frac{2 N p}{N + \alpha}} (\Rset^N)\) of the Choquard equation \eqref{eqAutoHomoChoquard} such that \(u \in W^{2, 2}_\mathrm{loc}(\Rset^N) \cap W^{1,\frac{2 N}{N + \alpha} p} (\Rset^N)\)
satisfies the Poho{\v{z}}aev identity \citelist{\cite{Menzala1983}*{(2.8)}\cite{CingolaniSecchiSquassina2010}*{lemma 2.1}\cite{MorozVanSchaftingen2013JFA}*{proposition~3.1}\cite{MorozVanSchaftingen2015TAMS}*{theorem~3}\cite{GenevVenkov2012}*{(56)}}
\begin{equation}\label{e-Pohozaev}
  \frac{N - 2}{2} \int_{\Rset^N} \abs{\nabla u}^2 + \frac{N}{2} \int_{\Rset^N} \abs{u}^2
= \frac{N + \alpha}{2 p} \int_{\Rset^N} \bigl(I_\alpha \ast \abs{u}^p\bigr) \abs{u}^{p}.
\end{equation}
The proof goes as in the local case \citelist{\cite{Pohozaev1965}\cite{PucciSerrin1986}} by testing the equation against suitable cut-offs of the function
\(x \in \Rset^N \mapsto x \cdot \nabla u (x) \in \Rset\).

The global integrability and local regularity assumptions are satisfied by any weak solution
in the space \(H^1 (\Rset^N)\) as soon
as \(\frac{N - 2}{N + \alpha} \le \frac{1}{p} \le \frac{N}{N + \alpha}\) \cite{MorozVanSchaftingen2015TAMS}*{theorem~2}.

\subsubsection{Numerical computations}

In the case \(N = 3\), \(\alpha = 2\) and \(p = 2\), the energy levels have been computed numerically for groundstates \cite{RuffiniBonazzola1969} and for radial boundstates \citelist{\cite{BernsteinGidaliJones1998}\cite{MorozPenroseTod1998}}.
The relationship of these numerical computations with formal computations \cite{KumarSoni2000} was not clear; the numerical computations have been supported by new bounds \cite{Tod2001}.

\subsection{Further equations}

\subsubsection{Nonhomogeneous potentials}
The existence of solutions has been studied for Choquard equations \eqref{eqAutoHomoChoquard} with \(p = 2\) when
the Riesz potential is replaced by \emph{a more general interaction potential} \(W : \Rset^N \to [0, \infty)\):
\begin{equation}
 \label{eqAutoHomoChoquard-W}
 -\Delta u + u = \bigl(W \ast \abs{u}^2\bigr) u \qquad \text{in \(\Rset^N\)}.
\end{equation}

In the case where \(N = 3\) and
\(W\) is a \emph{Yukawa potential} defined for each \(x \in \Rset^N \setminus \{0\}\)
by \(W (x) = \exp(-\abs{x})/(4 \pi \abs{x})\), which
is the fundamental solution of the linear operator \(-\Delta + 1\) on \(\Rset^3\) and which is
a special case of the \emph{Bessel potential}, the existence of solutions was proved
by ordinary differential equations methods \cite{KarasevMaslov1979}.
More generally if \(N \le 5\) and $W=B_\lambda$ is the modified Bessel potential,
that is the Green function of the Helmholtz operator $(-\Delta+\lambda I)^{-1}$ with $\lambda>0$,
there exists a groundstate and infinitely many radial solutions \cite{ZhaoZhaoShi2015}; the assumption \(N \ge 5\) is sharp in view of a Poho\v zaev identity \cite{ZhaoZhaoShi2015}*{(2.2)}. For $3\le N\le 5$ and $\lambda\to 0$, the groundstates of \eqref{eqAutoHomoChoquard-W}
converge (up to a subsequence) to the groundstate of the Choquard equation \eqref{eqAutoHomoChoquard} with $\alpha=2$ and $p=2$,
while for $1\le N\le 3$ and $\lambda\to\infty$ the groundstates of \eqref{eqAutoHomoChoquard-W} converge to the groundstates
of the local equation $-\Delta u+u=u^3$ in $\Rset^N$ \cite{ZhaoZhaoShi2015}*{theorem~1.3}.

When \(W \in L^{q_1} (\Rset^N) + L^{q_2} (\Rset^N)\) and \(W \ge 0\) solutions were constructed
by various variational methods \citelist{\cite{Lions1980}\cite{Lions1982}\cite{Menzala1980}}.

\subsubsection{Low dimensional Choquard equation}
The case where the interaction potential \(W\) is the \emph{Newtonian potential} in the one- or two-dimensional space
is particularly interesting. Indeed one has then for \(x \in \Rset\), \(W (x) = c - \abs{x}/2\)
or for \(x \in \Rset^2\), \(W (x) = (c - \log \abs{x})/(2 \pi)\):
the potential \(W\) changes sign and moreover its negative part is unbounded.
In particular, the action functional is not anymore well-defined on the Sobolev space
\(H^1 (\Rset^N)\). An idea could be to take the quantity
\begin{equation}
\label{eqNegPotNorm}
 \Bigl(\int_{\Rset^N} \abs{\nabla u}^2 + \abs{u}^2\Bigr)^\frac{1}{2}
 + \Bigl(\int_{\Rset^N} \bigl(W \ast \abs{u}^2\bigr) \abs{u}^2\Bigr)^\frac{1}{4}
\end{equation}
as a norm \cite{Lions1980}. Such ideas have been implemented successfully for a class of
Schr\"odinger--Poisson--Slater equation \cite{Ruiz2010}
(see also \cite{MercuriMorozVanSchaftingen}).
When \(W_-\) is not bounded, it can be observed that the quantity \eqref{eqNegPotNorm} \emph{is not a norm}:
large enough translations of a fixed compactly supported smooth function do not satisfy the triangle inequality.

By considering the subspace on which \(\int_{\Rset^N} W_- \abs{u}^2 < \infty\),
\emph{the existence and uniqueness of groundstates and the existence of bound states} has been proved when \(W\) is a low-dimensional
Newtonian potential \citelist{\cite{ChoquardStubbe2007}
\cite{Stubbe2008}\cite{StubbeVuffray2010}}. In order to handle the absence of invariance under translation
of the norm on this restricted natural space, the problem is reduced to the radial case
by symmetrisation arguments.
Solutions can also be constructed for such problems by ordinary differential equation
methods \cite{ChoquardStubbeVuffray2008}.

High energy solutions have also been constructed variationally \cite{CingolaniWeth}.
The construction required to manage the combination of a translation invariant functional
with a non-translation-invariant norm.

\subsubsection{Constant magnetic field}
\label{sectionConstantMagneticField}
When \(A : \Rset^N \to \Rset^N\) is a skew-symmetric linear map, the \emph{magnetic Choquard equation}
prescribes for \(u : \Rset^N \to \Cset\)
\begin{equation}
\label{eqMagAutoChoquard}
 (-i\nabla + A)^2 u + u = \bigl(I_\alpha \ast \abs{u}^p\bigr) \abs{u}^{p - 2} u.
\end{equation}
This problem is invariant under the noncommutative group of phase rotations defined for \(\alpha \in \Cset\) and \(\abs{\alpha} = 1\) by the action \(u
\mapsto \alpha u\) and magnetic translations defined for \(h \in \Rset^N\)
by the action
\[
  u \longmapsto \bigl(x \mapsto u (x - h) e^{-i A(h)\cdot (x-h/2)}\bigr).
\]
In particular the action functional \(\mathcal{A}\) corresponding to the magnetic Choquard equation \eqref{eqMagAutoChoquard} is invariant under a noncompact group which is locally compact.

The existence of a groundstate has been proved in the three-dimensional case \(N = 3\), \(\alpha = 2\) and \(p = 2\) \cite{GriesemerHantschWellig2012}*{theorem~2.3}.
If \(\frac{N - 2}{N + \alpha} < \frac{1}{p} < \frac{N}{N + \alpha}\)
and if \(\dim \ker A \ne 1\), then there exists infinitely many solutions \cite{CingolaniClappSecchi2012}*{theorem~1.1}.
The assumption on the magnetic potential \(A\) is satisfied
if either the dimension \(N\) is even or if \(\dim \ker A \ge 3\).
(In particular no nontrivial constant magnetic is covered been covered in the
three-dimensional case \(N = 3\).)

The asymptotics of the groundstate level have been studied when \(N = 3\), \(\alpha = 2\) and \(p = 2\)
when the magnetic field \(\operatorname{curl} A\) becomes large \cite{FrankGeisinger2015}*{theorem~1.2}.

\subsubsection{General nonlinearity}
\label{sectionGeneralNonlinearity}
The Choquard equation has also been studied
when the homogeneous nonlinearities \(\abs{u}^p\) and \(\abs{u}^{p - 2} u\) is replaced by a \emph{more general nonlinearity}. If one wants to keep the variational structure,
the equations writes then as
\begin{equation}
\label{problemChoquardBL}
 - \Delta u + u = \bigl(I_\alpha \ast F (u)\bigr) F' (u) \quad \text{in \(\Rset^N\)}.
\end{equation}
When \(N \ge 3\), if for each \(s \in \Rset\),
\[\abs{s F' (s)} \le C \bigl(\abs{s}^\frac{N + \alpha}{N} + \abs{s}^\frac{N + \alpha}{N - 2}\bigr),\]
\(\lim_{s \to 0} {F (s)}/{\abs{s}^\frac{N + \alpha}{N}} = 0\), \(\lim_{\abs{s} \to \infty} {F (s)}/{\abs{s}^\frac{N + \alpha}{N - 2}} = 0,\) and if there exists \(s_0 \in \Rset \setminus \{0\}\) such that \(F (s_0) \ne 0\), then the general nonlinear Choquard equation \eqref{problemChoquardBL} has a solution \cite{MorozVanSchaftingen2015TAMS}.
A similar results holds in the planar case \(N = 2\) with possible exponential growth of the nonlinearity \cite{BattagliaVanSchaftingen}. Solutions of the corresponding normalised problem, that is, when the total mass \(\int_{\Rset^N} \abs{u}^2\) is prescribed and a Lagrange multiplier is allowed have also been constructed \cite{LiYe2014}.

This result is a counterpart for the Choquard equation of the classical
result of Berestycki and Lions \cite{BerestyckiLions1983} (see also \cite{BerestyckiGallouetKavian1983}), the assumptions are similarly
almost necessary. The proof uses a variational trick related to the Poho\v zaev identity
\citelist{\cite{Jeanjean1997}\cite{HirataIkomaTanaka}}
and does not
seem to work directly for a more general potential \(W\) in place of the Riesz potential
\(I_\alpha\).

Existence of solutions has also been proved for exponential nonlinearities, under an Ambrosetti--Rabinowitz superlinearity assumption \cite{AlvesCassaniTarsiYang}.

\subsubsection{Local autonomous perturbation}%
\label{sectLocalQ}

The Choquard equation \eqref{eqAutoHomoChoquard} can be perturbed by a local nonlinearity
\[
 -\Delta u + u = \bigl(I_\alpha \ast \abs{u}^p\bigr) \abs{u}^{p - 2} u + \abs{u}^{q - 2} u\qquad \text{in \(\Rset^N\)}.
\]
The existence of solutions has been proved when \(N = 3\), \(0 < \alpha < 1\), \(p = 2\)
and \(4 \le q < 6\) \cite{ChenGuo2007}.
When \(N = 2\), the Riesz potential \(I_\alpha\) is replaced by the two-dimensional Newtonian potential \(W\) defined for \(x \in \Rset^2 \setminus \{0\}\) by \(W (x) = -\log (\abs{x})/(2 \pi)\),
then the problem has been studied for \(p = 2\) and \(q > 2\) \cite{CingolaniWeth}.

\subsubsection{Nonvariational case}
\emph{Self-dual variational principles} have allowed to treat the nonvariational problem
\begin{equation}
 \label{eqNonVarChoquard}
 -\Delta u + u = \bigl(W \ast \abs{u}^p\bigr) \abs{u}^{q - 2} u \qquad \text{in \(\Rset^N\)},
\end{equation}
when \(W \in L^1 (\Rset)^N\), \(1 \le p < \frac{N}{N - 2}\), \(1 < q < \frac{N}{N - 2}\) and \(pq < 2\) \cite{Ghoussoub2009}*{theorem~12.5}.

\subsubsection{Pseudorelativistic and fractional Choquard equation}
The pseudorelativistic Choquard equation
\begin{equation}
 \label{eqPseudoChoquard}
 \sqrt{-\Delta+m^2\,}\; u + \mu u = \bigl(I_2 \ast \abs{u}^2\bigr)u \qquad \text{in \(\Rset^3\)}
\end{equation}
appears as a model of pseudo-relativistic boson stars  in the mean-field limit \citelist{\cite{LiebYau1987}\cite{FrohlichJonssonLenzmann2007}}.
This equation and other fractional modifications of the Cho\-quard equation have been studied in \citelist{\cite{Lenzmann2009}\cite{ZelatiNolasco11}\cite{Mugnai13}\cite{CingolaniSecchi2015}\cite{DAveniaSquassina2015}},
see also further references therein.

\section{Nonautonomous equation}
\label{sectionNonautonomous}

In this section we explore Choquard equations which are not invariant under translations of the Euclidean space \(\Rset^N\), due to the presence of a variable electric potential or of a magnetic potential in the linear part of the equation. This can also happen through a non translation-invariant modifications of the nonlinear term.

\subsection{Electric potential}

Physical models in which particles are under the influence of an \emph{external electric field},
lead to study Choquard equations in the form
\begin{equation}
 \label{eqChoquard-V}
 -\Delta u + Vu = \bigl(I_\alpha \ast \abs{u}^p\bigr) \abs{u}^{p - 2} u \qquad \text{in \(\Rset^N\)},
\end{equation}
where $V\in L^1_{\mathrm{loc}}(\Rset^N)$ is a non-constant electric potential.

\subsubsection{Perturbation of a constant potential}
P.--L.~Lions \cite{Lions1984-1}*{\S 3} has studied the existence of solutions to problem \eqref{eqChoquard-V} with $p=2$ when the external potential $V$
is a \emph{perturbation of a constant potential}, that is, for each \(x \in \Rset^N\), \(V\) can be written as
\[
  V(x)=1+V_0(x),
\]
where $V_0$ decays at infinity and satisfies some mild regularity assumptions, for example
$V_0\in L^\frac{N}{2} (\Rset^N)+L^t(\Rset^N)$ with $\frac{N}{2}\le t<\infty$ if $N\ge 3$.
To construct solutions of the equation \eqref{eqChoquard-V}, Lions has considered a maximisation problem
$$
  \mathcal I_V:=\sup\Big\{\int_{\Rset^N} \bigl(I_\alpha \ast \abs{u}^2\bigr)\abs{u}^2\st u \in H^1 (\Rset^N),\int_{\Rset^N} \Big(\abs{\nabla u}^2+V\abs{u}^2\Big)=1\Big\}.
$$
Up to multiplication by a constant, maximisers of $\mathcal I$ are solutions of Choquard equation \eqref{eqChoquard-V}.
The associated \emph{limit problem at infinity}
$$
\mathcal I^\infty:=\sup\Big\{\int_{\Rset^N} \bigl(I_\alpha \ast \abs{u}^2\bigr)\abs{u}^2\st u \in H^1 (\Rset^N),\int_{\Rset^N} \Big(\abs{\nabla u}^2+\abs{u}^2\Big)=1\Big\},$$
is merely a reformulation of the minimisation of the Sobolev quotient \eqref{e-min-Sobolev} associated to the Choquard equation, which for $p=2$ admits a minimiser if and only if
$(N-4)_+<\alpha<N$.

As an applications of his concentration-compactness method, Lions has proved \cite{Lions1984-1}*{theorem~III.3} that
if $V_0\not\equiv 0$ and if there exists $\nu>0$ such that
$$\int_{\Rset^N} \bigl(\abs{\nabla u}^2+V\abs{u}^2\bigr)\ge \nu\|u\|_{H^1(\Rset^N)}^2\qquad \text{for each \(u\in H^1(\Rset^N)\)}.$$
then every maximising sequence for $\mathcal I$ is relatively compact in $H^1(\Rset^N)$ if and only if
$$
  \mathcal I_V<\mathcal I^\infty.
$$
Similar results, as well as results on the existence of $L^2$-constrained minimisers, have been obtained for the Choquard equations
\begin{equation}%
\label{eqChoquard-VW}%
-\Delta u+Vu= \bigl(W*\abs{u}^2\bigr)u\qquad\text{in $\Rset^N$},
\end{equation}
with a general, possibly sign--changing, convolution kernels $W$ \cite{Lions1984-1}*{theorems III.1 and III.3}.
Unlike E.\thinspace{}Lieb's approach \cite{Lieb1977}, which relies on the symmetrisation techniques, P.-L.\thinspace{}Lions's results do not require any symmetry properties of the potentials $V$ and $W$.

Although P.-L.\thinspace{}Lions has stated his results \cite{Lions1984-1} for $p=2$ and $(N-4)_+<\alpha<N$ only,
they can be extended to the same noncritical range as in theorem~\ref{theoremExistenceGroundStates}
$p\in\big(\frac{N+\alpha}{N},\frac{N+\alpha}{(N-2)_+}\big)$.
The endpoints of this interval require separate consideration.

\subsubsection{Lower and upper critical exponents}
We have seen above in proposition~\ref{propositionAutoHomogFunctional} that in the case of a constant potential $V(x)\equiv 1$,
the action functional $\mathcal A$ is a well-defined continuously Fr\'echet differentiable functional on the Sobolev space $H^1(\Rset^N)$ if and only if $p\in\big[\frac{N+\alpha}{N},\frac{N+\alpha}{(N-2)_+}\big]$.
However, as a consequence of the Poho{\v{z}}aev identity \eqref{e-Pohozaev}, the autonomous Choquard equation \eqref{eqAutoHomoChoquard}
with critical exponents $p=\frac{N+\alpha}{N-2}$ and $p=\frac{N+\alpha}{N}$ does not have any nontrivial solution in the Sobolev space $H^1(\Rset^N)$.
In these two critical cases, it would be meaningless to study the nonautonomous equation \eqref{eqChoquard-V} as a perturbation of the autonomous equation \eqref{eqAutoHomoChoquard}.

The Choquard equation with the \emph{lower critical exponent},
\begin{equation}
\label{equationNLChoquard}
 - \Delta u + V u  = \bigl(I_\alpha \ast \abs{u}^{\frac{\alpha}{N}+1}\bigr) \abs{u}^{\frac{\alpha}{N}-1} u\qquad\text{in \(\Rset^N\)},
\end{equation}
has been studied by the authors \cite{MorozVanSchaftingen2015CCM}.
The exponent $\frac{\alpha}{N}+1$ is critical with respect to the Hardy--Littlewood--Sobolev inequality \eqref{eqHLS}, which can be reformulated variationally as
\[
\mathcal J_\infty =  \inf \Bigl\{ \int_{\Rset^N} \abs{u}^2 \st u \in L^2 (\Rset^N), \int_{\Rset^N} \bigl(I_\alpha \ast \abs{u}^{\frac{\alpha}{N} + 1}\bigr) \abs{u}^{\frac{\alpha}{N} + 1} = 1\Bigr\} > 0.
\]
This infimum \(\mathcal J_\infty\) is achieved if and only if
\begin{equation}\label{HLSoptimal}
  u(x) = C\left(\frac{\lambda}{\lambda^2 + \abs{x - a}^2}\right)^{N/2},
\end{equation}
where \(C > 0\) is a fixed constant, \(a \in \Rset^N\) and \(\lambda \in (0, \infty)\) are parameters \citelist{\cite{Lieb1983}*{theorem~3.1}\cite{LiebLoss2001}*{theorem~4.3}}.
The presence of this nonlinear lower critical exponent is a feature of the Choquard equation that does not appear in its local counterpart the nonlinear Schr\"odinger equation.

The existence of nontrivial solutions for \eqref{equationNLChoquard}, can be obtained by considering the minimisation problem
\[
  \mathcal J_V = \inf \Bigl\{ \int_{\Rset^N} \bigl(\abs{\nabla u}^2 + V \abs{u}^2\bigr) \st u \in H^1 (\Rset^N), \int_{\Rset^N} \bigl(I_\alpha \ast \abs{u}^{\frac{\alpha}{N} + 1}\bigr) \abs{u}^{\frac{\alpha}{N} + 1} = 1\Bigr\}.
\]
Up to multiplication by a constant, minimisers of $\mathcal J_V$ are solutions of the Choquard equation \eqref{equationNLChoquard}.
Direct substitution of minimisers of the form \eqref{HLSoptimal} with \(\lambda\to \infty\) into $\mathcal J_V$ shows that if $V\equiv 1$, then
$$\mathcal J_V=\mathcal J_\infty,$$
so that $\mathcal J_\infty$ could indeed be seen as a limit problem at infinity for $\mathcal J_V$.
The form of minimisers in \eqref{HLSoptimal} suggests that a loss of compactness in $\mathcal J_V$ may occur by translations and dilations at infinity.
By a Brezis-Lieb type lemma for Riesz potentials \eqref{eqBrezisLieb}
and a concentration-compactness argument, it has been showed that if $V\in L^\infty(\Rset^N)$ and $\lim_{|x|\to\infty} V(x)= 1$ then
every minimising sequence for $\mathcal J_V$ is relatively compact in $H^1(\Rset^N)$ if and only if \cite{MorozVanSchaftingen2015CCM}*{theorem~3 and proposition 5}
\begin{equation}\label{e-minimising-J}
\mathcal J_V<\mathcal J_\infty.
\end{equation}
Moreover, if
\begin{equation}
\liminf_{\abs{x} \to \infty} \bigl(1 - V (x)\bigr)\abs{x}^2 > \frac{N^2 (N - 2)}{4 (N + 1)}
\end{equation}
then $\mathcal J_V<\mathcal J_\infty$ and hence Choquard equation \eqref{eqChoquard-V} has a nontrivial positive solution.
Some necessary conditions on the potential $V$ for the strict inequality to take place \eqref{e-minimising-J} have been discussed.

The Choquard equation with the \emph{upper critical exponent} is
\begin{equation}
\label{equationNLChoquard-uppercritical}
 - \Delta u + V u  = \bigl(I_\alpha \ast \abs{u}^{\frac{N+\alpha}{N-2}}\bigr) \abs{u}^{\frac{\alpha - N + 4}{N-2}} u \qquad\text{in \(\Rset^N\)}.
\end{equation}
The exponent $\frac{N+\alpha}{N-2}$ is critical with respect to the Hardy--Littlewood--Sobolev inequality in the form
\[
\mathcal J^\infty =  \inf \Bigl\{ \int_{\Rset^N} \abs{\nabla u}^2 \st u \in D^1(\Rset^N), \int_{\Rset^N} \bigl(I_\alpha \ast \abs{u}^\frac{N+\alpha}{N-2}\bigr) \abs{u}^\frac{N+\alpha}{N-2} = 1\Bigr\} > 0.
\]
Minimisers for \(\mathcal J^\infty\)
are known explicitly \citelist{\cite{Lieb1983}*{theorem~3.1}\cite{LiebLoss2001}*{theorem~4.3}}].
Direct substitution of minimisers for $\mathcal J^\infty$ into $\mathcal J_V$ shows that
$\mathcal J_V=\mathcal J^\infty$, provided that the potential $V$ is positive on an open subset of $\Rset^N$,
so the condition $\mathcal J_V<\mathcal J^\infty$ is no longer meaningful,
and \eqref{equationNLChoquard-uppercritical} is a Brezis--Nirenberg \cite{BrezisNirenberg1983} type problem.

Although \eqref{equationNLChoquard-uppercritical} has not yet been studied, existence results for \eqref{equationNLChoquard-uppercritical} on bounded domains $\Omega\subset\Rset^N$ have been recently obtained \citelist{\cite{GaoYang2016}\cite{MukherjeeSreenadh2016}}, as well as some perturbations by a local nonlinear term \cite{GaoYang2016b}.
We bring to the attention of the reader that other versions of the Choquard equation can be defined on a bounded domain. For example, the Riesz potential \(I_\alpha\) could be replaced by \((-\Delta)^{-\alpha/2}\), where \(\Delta\) is the Laplacian on \(\Omega\) with Dirichlet boundary conditions.
The physical and mathematical relevance of the various possible integral kernels on a domain deserves some study in the future.

\begin{problem}
Under which conditions on the potential $V$ does \eqref{equationNLChoquard-uppercritical} have a positive solution?
\end{problem}

Starting points could be the methods used to treat critical problems on the Euclidean space in \(\Rset^N\) by either variational concentration-compactness methods \cite{BenciCerami1990} or perturbation methods \cite{Ambrosettietal1999}.

\subsubsection{Confining potentials}\label{sectConfining}
The action functional \(\mathcal{A}_V\) for equation \eqref{eqChoquard-V} is defined for \(u : \Rset^N \to \Rset\) by
\begin{equation}
\label{eqAutoHomoAction-V}
 \mathcal{A}_V (u)
 = \frac{1}{2} \int_{\Rset^N} \bigl(\abs{\nabla u}^2 + V\abs{u}^2\bigr)
 - \frac{1}{2 p} \int_{\Rset^N} \bigl(I_\alpha \ast \abs{u}^{p}\bigr) \abs{u}^p.
\end{equation}
If we assume that $V$ is a nonnegative \emph{confining} external potential, that is, if
$$\lim_{|x|\to\infty}V(x)=+\infty,$$
then the natural domain for the functional $\mathcal{A}_V$ is the weighted Sobolev space $H^1_V(\Rset^N)$,
which is the completion of the class $C^\infty_c(\Rset^N)$ with respect to the norm
$\|u\|_V:=\big(\int_{\Rset^N} \abs{\nabla u}^2 + V\abs{u}^2\big)^{1/2}$.

In the case of a confining potential \(V\), the space $H^1_V(\Rset^N)$ is compactly embedded into $L^2(\Rset^N)$;
this simplifies considerably the analysis of the Palais--Smale sequences for $\mathcal A_V$.
The existence of a positive radially symmetric groundstate solution for \eqref{eqChoquard-V} with a \emph{radial} confining potential
$V$ has been studied when \(p = 2\) \cite{FrohlichLenzmann2004}. When \(\alpha = 2\), see also \cite{CaoWangZou2012}.
The specific case of harmonic potential $V(x)=|x|^2$ and $p\in\big(\frac{N+\alpha}{N},\frac{N+\alpha}{N-2}\big)$ has been studied in \cite{Feng2016}. (Explicit formulas for solutions with a quadratic external potential with a quadratic nonlocal potential (corresponding formally to \(\alpha = -2\)) have been given \cite{MaedaMasaki2013}.)

Strongly enough confining potentials even allow to \emph{enlarge the admissibility range} for the nonlinearity exponent \(p\).
Although the well--posedness interval of proposition~\ref{propositionAutoHomogFunctional} is no longer valid for $\mathcal A_V$, the Stein--Weiss weighted Hardy--Littlewood--Sobolev type inequality \eqref{eqHLSSW} has been used to prove for example, that, if \(V\) is confining  \(p \in (1, \frac{N+\alpha}{N}]\) and if
\begin{equation}
\label{eqVSX}
 \lim_{|x|\to\infty} \frac{V(x)}{|x|^{\frac{N + \alpha}{p} - N}} = \infty,
\end{equation}
then the functional $\mathcal A_V$ is well-defined and continuously Fr\'echet--differentiable on $H^1_V(\Rset^N)$.
This has allowed to prove existence of groundstates if either \( \frac{N + \alpha}{N} < p \frac{N + \alpha}{(N - 2)_+}\) and \(V\) is confining, or \(1 <p < \frac{N+\alpha}{N}\) and \eqref{eqVSX} holds \cite{VanSchaftingenXia}.
This result is a counterpart of the classical result for the nonlinear Schr\"odinger equation \cite{Rabinowitz1992}*{theorem~1.7}.

\subsubsection{Periodic potentials}
The local compactness of Pa\-lais--Smale sequences can be also obtained relatively easy
when the potential \(V\) is periodic.
The existence of a nontrivial weak solution for Choquard equations
\begin{equation}\label{eqChoquard-VWf}
-\Delta u+Vu=\bigl(W*F(u)\bigr)F'(u)\qquad\text{in $\Rset^N$},
\end{equation}
with a \emph{positive periodic electric potential} $V$ and under a general assumptions on the even convolution kernel $W\ge 0$ and nonlinearity $F$
has been proved by a Mountain Pass type argument \cite{Ackermann2004}*{theorem~2.1}.
If, in addition, the function $F$ is even,
then there exists infinitely many geometrically distinct solutions.

The situation is more involved when the periodic potential $V$ is allowed to change sign.
In this case the operator $-\Delta+V$ may have some essential spectrum below $0$.
As a consequence, the quadratic part of the action functional for \eqref{eqChoquard-VWf} is \emph{strongly indefinite}
and the analysis of Palais--Smale sequences becomes much more delicate.
The existence of at least one nontrivial solution for \eqref{eqChoquard-V} with a sign changing periodic $V$ in the case $N=3$, $\alpha=2$, $p=2$ has been proved under the assumption that $0$ is in the gap of the spectrum of the operator $-\Delta+V$
\cite{BuffoniJeanjeanStuart1993}.
The result has been extended to solutions of the problem \eqref{eqChoquard-VWf} and multiplicity results have been established when the nonlinearity \(F\) is an even function \cite{Ackermann2004}*{theorem~2.2}.

Choquard equations with \emph{mixed} potentials $V$ which are periodic with respect to some dimensions and confining with respect to the others have been considered \cite{SoutodeLima}.

In low dimensions, the existence of solutions to the \emph{planar logarithmic Choquard equation}
\begin{equation*}
-\Delta u+Vu+\bigl(\log|x|*|u|^2\bigr)u=c|u|^{p-2}u\qquad\text{in $\Rset^2$},
\end{equation*}
under a positive periodic potential $V$ and with $c\ge 0$ has been obtained \cite{CingolaniWeth}.

In the one-dimensional case, \emph{periodic solutions} have been constructed for
for the periodic problem in \(\Rset\) \cite{Dymarksii1997}.

\subsubsection{Decaying potentials}
If $V\equiv 0$ or if $\lim_{|x|\to\infty}V(x)=0$ then $0$ belongs to the essential spectrum for $-\Delta+V$.
The space $H_V(\Rset^N)$ can be formally defined as before, although some care is needed in lower dimensions:  if $V\equiv 0$ and $N=1,2$ then $H^1_V(\Rset^N)$ is not continuously embedded into the space of distributions $\mathcal D^\prime(\Rset^N)$ (see for example \cite{PinchoverTintarev2006}*{p.66}).
A Hardy type inequality \cite{MorozVanSchaftingen2010}*{lemma 6.1} ensures that the space $H^1_V(\Rset^N)$ is well defined and continuously embedded into $L^2_{\mathrm{loc}}(\Rset^N)$
in dimensions $N=1,2$, provided that $V(x)>0$ on an open subset of $\Rset^N$.
However, if $V\equiv 0$ or if \(V\) compactly supported, the action functional $\mathcal A_V$ is typically not well-defined nor Fr\'eschet differentiable on $H_V(\Rset^N)$.
This difficulty is not only technical. Liouville type theorems show that Choquard equations with fast decaying potentials
indeed do not have positive solutions or even positive super-solutions for certain ranges of parameters.

\subsubsection{Liouville theorems}%
\label{sLiouville}

Nonlinear Liouville theorems state the \emph{nonexistence of positive solutions}
of elliptic equations or inequalities.
A typical nonlinear Liouville theorem in the local case says that the inequality
$$
  -\Delta u\ge u^p\qquad \text{in \(\Rset^N\setminus B_\rho\)}
$$
admits a positive classical solution if and only if $p>\frac{N}{N-2}$. This result goes back to J. Serrin in the 1970s (see \cite{Quittner-Souplet-2007}*{theorem~8.4} for a proof and references).
The value of the critical exponent $\frac{N}{N-2}$ does not depend on the value of $\rho>0$. The statement is robust:
it still holds if we perturb $-\Delta$ by a sufficiently weak potential $V$,
for example a potential that satisfies the bound $|V(x)|\le c|x|^{-2-\delta}$, for some $\delta>0$.

Liouville type theorems for the Choquard inequalities
\begin{equation}\label{e-Liouville-V}
-\Delta u+Vu\ge \bigl(I_\alpha\ast u^p\bigr)u^q\qquad \text{ in }\Rset^N\setminus B_\rho,
\end{equation} have been obtained \cite{MorozVanSchaftingen2013JDE}.
We assume in \eqref{e-Liouville-V} that $p>0$ and $q\in\Rset$,
and  we observe that equation \eqref{e-Liouville-V-0} has no variational structure, unless $q=p-1$.
By $I_\alpha\ast f$ in a domain $\Omega\subset\Rset^N$ we understand
the convolution $I_\alpha*(\chi_\Omega f)$ in $\Rset^N$, restricted to $\Omega$.
Note that the concept of a solution to \eqref{e-Liouville-V} in a domain is nonlocal.
In particular, if $u\ge 0$ is a solution of Choquard equation
\begin{equation}\label{e-Liouville-V-0}
-\Delta u+Vu= \bigl(I_\alpha\ast u^p\bigr)u^q\qquad \text{in }\Rset^N,
\end{equation}
then $u$ is only a \emph{super}-solution
of the same equation in a proper subdomain $\Omega\subsetneq\Rset^N$.
From this point of view, consideration of \emph{inequalities} \eqref{e-Liouville-V} is quite natural.
Generally speaking, it is desirable to obtain nonexistence results for as wide as  possible class of solutions.
Throughout this section by a solution of \eqref{e-Liouville-V} we understand a nonnegative function $u\in L^p(\Rset^N\setminus B_\rho,|x|^{-(N-\alpha)}dx)\cap L^1_{\mathrm{loc}}(\Rset^N\setminus B_\rho,V(x)\dif x)$ which satisfies \eqref{e-Liouville-V} in the distributional sense.
Here we must assume that $u\in L^p(\Rset^N\setminus B_\rho,|x|^{-(N-\alpha)}dx)$ in order to have $I_\alpha\ast u^p<+\infty$ almost everywhere in $\Rset^N\setminus B_1$, see \eqref{e-Riesz-necessary}.

The nonexistence arguments were based on a \emph{nonlocal positivity principle} \cite{MorozVanSchaftingen2013JDE}*{proposition~3.2},
which claims that if $u>0$ is a distributional solution of \eqref{e-Liouville-V} with arbitrary, possibly sign--changing potentials $V\in L^1_{\mathrm{loc}}(\Rset^N)$,
then \(u^{q - 1} \in L^1_{\mathrm{loc}} (\Rset^N\setminus B_\rho)\) and
for all \(R>\rho\) and for all \(\varphi \in C^\infty_c(B_{4R}\setminus B_\rho)\) the following inequality holds
\begin{equation}\label{e-AAP}
\int_{\Rset^N\setminus B_\rho}\big(|\nabla \varphi|^2+V(x)\varphi^2\big)\ge
\frac{C_\alpha}{R^{N - \alpha}} \Bigl(\int_{B_{4R}\setminus B_\rho} u^{p}\Bigr) \Bigl(\int_{\Rset^N\setminus B_\rho} u^{q - 1} \varphi^2 \Bigr).
\end{equation}
This statement shall be seen as a nonlocal version of the Agmon--Allegretto--Piepenbrink positivity principle \cite{CyconFroeseKirschSimon}*{theorem~2.12}.
Indeed, if $u$ is a positive solution of \eqref{e-Liouville-V} then \eqref{e-AAP} claims that the quadratic form associated to the Schr\"odinger operator $-\Delta+V$ in $\Rset^N\setminus B_\rho$ is nonnegative definite and admits a lower bound (in the right hand side of \eqref{e-AAP}),
which encodes the information about the order of the Riesz potential $I_\alpha$ and the decay of the function $u$.

To give a taste of the techniques of nonexistence proofs, consider the case of the autonomous Choquard inequality
\begin{equation}\label{e-Liouville}
-\Delta u\ge\bigl(I_\alpha \ast u^p\bigr)u^q\qquad \text{in }\Rset^N\setminus B_\rho.
\end{equation}
Assume that $u\ge 0$ is a solution of \eqref{e-Liouville}.
First of all, note that then $u$ is a superharmonic function, that is, $-\Delta u\ge 0$ in $\Rset^N\setminus B_\rho$.
If $u\not\equiv 0$, by comparison with the harmonic function $|x|^{-(N-2)}$
we obtain a lower bound
\begin{equation}\label{e-harmonic}
u(x) \ge c|x|^{-(N-2)},\qquad x\in\Rset^N\setminus B_{2\rho}.
\end{equation}
If $p\le \frac{\alpha}{N-2}$, \eqref{e-harmonic} means that $u\not\in L^p(\Rset^N,|x|^{-(N-\alpha)}dx)$.
We conclude that $u\equiv 0$, or otherwise $I_\alpha*u^p\equiv +\infty$. This establishes the first \emph{nonexistence r\'egime}
for \eqref{e-Liouville}.

When $p> \frac{\alpha}{N-2}$ we will employ \eqref{e-AAP} with the family of test function $\varphi_R(|x|)=\varphi\big(\frac{|x|}{R}\big)$,
where $\varphi\in C^\infty_c(B_4\setminus\{0\})$ and
$\varphi_{|[1,2]}=1$. Then for all $R\gg\rho$ we deduce from \eqref{e-AAP} a so-called \emph{master inequality}
\begin{equation}\label{e-MI}
C R^{N -2}\ge R^{-(N - \alpha)}\Bigl(\int_{B_{2R}\setminus B_\rho} u^{p} \Bigr) \Bigl(\int_{B_{2R}\setminus B_R} u^{q - 1}\Bigr).
\end{equation}
Combined with the Cauchy--Schwarz inequality and with the lower bound \eqref{e-harmonic},
for $p+q\ge 1$ and all $R\gg\rho$ this leads to the estimate
\begin{multline*}
C R^{2N -2 - \alpha}\ge\Bigl(\int_{B_{2R}\setminus B_\rho} u^{p} \Bigr) \Bigl(\int_{B_{2R}\setminus B_R} u^{q - 1}\Bigr)\\
\ge \Bigl(\int_{B_{2R} \setminus B_R} u^\frac{p + q - 1}{2}\Bigr)^2 \ge c R^{2 N-(N - 2)(p + q - 1)}.
\end{multline*}
Comparing the left and right hand sides we conclude that $u\equiv 0$ provided that
$1\le p + q < \frac{N + \alpha}{N - 2}$.
Additional considerations allow to conclude that the critical line $p + q = \frac{N + \alpha}{N - 2}$
also belongs to the same nonexistence r\'egime, see \cite{MorozVanSchaftingen2013JDE}*{proposition~4.6}.

\begin{figure}[h!t]\label{f-L1}
\begin{center}
\includegraphics[width=.49\textwidth]{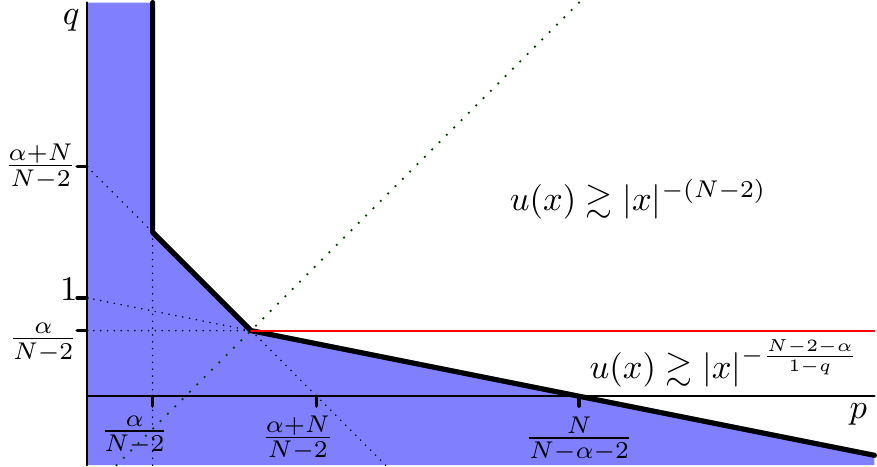}
\hspace{\stretch{1}}
\includegraphics[width=.49\textwidth]{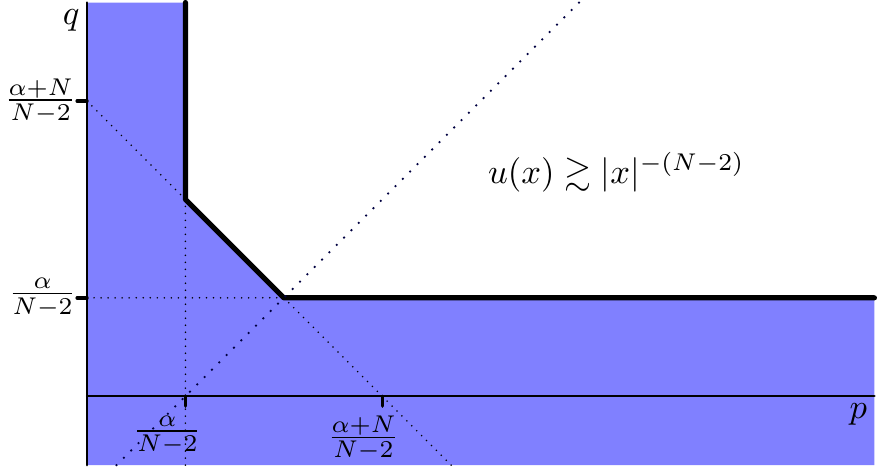}
\caption{Liouville map: region of existence of supersolution when $\alpha< N-2$ (left) and $\alpha\ge N-2$ (right), and maximal decay rates of supersolutions}
\end{center}
\end{figure}

The entire analysis of the $(p,q)$--plane for different sets of parameters is too long to reproduce here.
The final results are better represented by two \emph{Liouville maps}
of Choquard inequality \eqref{e-Liouville}, rather then by conventional theorem statements, see Figure~\ref{f-L1}.

%
%

We emphasise that \emph{all} the nonexistence r\'egimes on figure~\ref{f-L1} are obtained using only the master inequality \eqref{e-MI} and elementary bounds on positive superharmonic functions, although the use of \eqref{e-MI} varies significantly in different situations.
The Liouville maps on figure~\ref{f-L1} are sharp, which is confirmed by constructing explicit radial solutions for Choquard inequality \eqref{e-Liouville}, with the optimal decay rates indicated on
figure~\ref{f-L1}. For $q>\frac{\alpha}{N-2}$ the decay is controlled by the linear operator $-\Delta$.
For $q<\frac{\alpha}{N-2}$ and $\alpha<N-2$ the decay is governed by the nonlinear terms and becomes sensitive
to the exponent $q$. The decay bound $u\ge c|x|^{-\frac{N-2-\alpha}{1-q}}$
on figure~\ref{f-L1} is optimal and cannot be improved.
Optimal decay in the threshold case $q=\frac{\alpha}{N-2}$ is even more unusual and for $\alpha<N-2$ includes a logarithmic correction of the previous bound, see \cite{MorozVanSchaftingen2013JDE}*{proposition~4.11}.

For the readers' convenience we now summarise nonexistence results in the form of a conventional theorem,
but only in the variational case $q=p-1$, although it does not add any additional information to
Liouville maps of figure~\ref{f-L1}.

\begin{theorem}\label{Thm-free}
Let \(N\ge 3\), \(0<\alpha<N\), \(p > 0\), and \(\rho>0\).
Equation
\begin{equation}\label{e-Liouville-p}
-\Delta u\ge \bigl(I_\alpha\ast u^p\bigr)u^{p-1}\qquad \text{in }\Rset^N\setminus B_\rho
\end{equation}
has a nonnegative nontrivial solution in \(\Rset^N \setminus\Bar{B}_\rho\) if and only if:
\begin{align*}
p & > \frac{1}{2}\frac{2N - 2+\alpha}{N - 2} & & \text{if } 0<\alpha\le 2,\\
p & >  \frac{2N}{2N - \alpha - 2} & &\text{if } 2<\alpha< N - 2,\\
p & \ge 2 & &\text{if } 2<\alpha= N - 2,\\
p & > 1+\frac{\alpha}{N - 2} & &\text{if } \max\{2,N - 2\}<\alpha<N.
\end{align*}
\end{theorem}

So far we discussed only autonomous Choquard inequality \eqref{e-Liouville},
but Liouville results on Figure~\ref{f-L1} remain exactly the same for Choquard equations
\eqref{e-Liouville-V} with fast decaying potentials
$$V(x)=\frac{c}{|x|^\gamma},\qquad\gamma>2,$$
see \cite{MorozVanSchaftingen2013JDE}*{theorem~8}.
The same should remain true for any potential $V(x)$ with the property that
the Green functions of the operator $-\Delta+V$ is equivalent to $|x-y|^{-(N-2)}$.

For the Hardy type potentials
$$V(x)=\frac{\mu}{|x|^2}$$
with $\mu>-\frac{(N-2)^2}{4}$ the
Liouville properties remain qualitatively the same as on figure~\ref{f-L1},
but the values of critical exponents now depend explicitly on $\mu$,
see \cite{MorozVanSchaftingen2013JDE}*{theorem~9}.

In the case of slowly decaying potentials
$$V(x)=\frac{c}{|x|^\gamma},\qquad\text{with }\gamma<2,$$
which includes the constant potential $V(x)=c$, the results are different.
The primary reason for that is that the Green function of the operator $-\Delta+V$ decays exponentially,
which affects significantly all the previous arguments.
Here we only emphasise the fact that while for $p\ge 2$ the admissible decay of supersolutions to \eqref{e-Liouville}
with slowly decaying potentials is \emph{exponential}, for $p<2$ the decay rate is \emph{polynomial} as in \S\ref{s-Decay} above.
See \cite{MorozVanSchaftingen2013JDE}*{theorems 2-6} for details.

%
%

\subsubsection{Singular perturbations
and semiclassical limit}\label{sectSemiclass}

From the physical pro\-spective it is particularly important to study Choquard equations
\begin{equation}\label{e-V-eps}
-\eps^2\Delta u+V(x)u= \eps^{-\alpha}(I_\alpha\ast u^p)u^{p-1}\qquad \text{in }\Rset^N,
\end{equation}
where $\eps>0$ is a small parameter, typically related to the Planck constant.
Solutions of equation \eqref{e-V-eps} as $\eps\to 0$ are called \emph{semi-classical}.
Physically, it is expected that in the semi-classical limit $\eps\to 0$ there should be a correspondence between solutions of the equation \eqref{e-V-eps} and critical points of the potential $V$, which governs the \emph{classical} dynamics.

Mathematically, this can be justified by an observation that if $u_\eps$ is a solution of \eqref{e-V-eps} and $a\in\Rset^N$, then the
function $v_\eps(y)=u_\eps(a+\eps y)$ solves the rescaled equation
\begin{equation}\label{e-V-eps-rescaled}
-\Delta v_\eps+V(a+\eps y)v_\eps=(I_\alpha\ast v_\eps^p)v_\eps^q\qquad \text{in }\Rset^N.
\end{equation}
If $a \in \Rset^N$ is a critical point of the potential $V$ and $V(a)>0$ then the expectation is that $v_\eps$ should converge to a solution $v_0$ of the autonomous \emph{limit equation}
\begin{equation}\label{e-V-0}
-\Delta v_0+V(a)v_0=(I_\alpha\ast v_0^p)v_0^{p-1}\qquad \text{in }\Rset^N.
\end{equation}
If $v_0$ is a positive groundstate of \eqref{e-V-0}, constructed in theorem~\ref{theoremExistenceGroundStates},
then solution $u_\eps(x)\approx v_0\big(\frac{x-a}{\eps}\big)$ of the original equation should \emph{concentrate}
to $a$, in the sense that $u_\eps(x)\to 0$ as \(\eps \to 0\) for $x\neq a$ and $\liminf_{\eps \to 0} u_\eps(a)>0$.
For the local nonlinear Schr\"odinger equation
mathematical results of this type go back to Floer and Weinstein \cite{Floer-Weinstein} and by now well-understood.
An additional scaling parameter $\eps^{-\alpha}$, which appears in Choquard equation \eqref{e-V-eps} is required to ensure scaling invariance of the nonlocal problem.

First results on the existence of semiclassical solutions to Choquard equations of type \eqref{e-V-eps} have appeared
in the case $N=3$, $\alpha=2$ and $p=2$, under the assumption that $V\in C^2(\Rset^3)$ and $\inf_{x\in\Rset^3}V(x)>0$ \cite{WeiWinter2009}.
Using a \emph{Lyapunov--Schmidt reduction method}, the authors have proved that given non-degenerate critical points $a_1,\dots,a_m$ of the potential $V$,
there exist a family $v_\eps$ of multibump positive solutions which concentrate to that points.
Related results with a periodic external potential $V$ have also appeared \citelist{\cite{MacriNolasco2009}\cite{Nolasco2010}}. The existence of semiclassical solutions has been proved in the case $\liminf_{|x|\to\infty}V(x)|x|^\gamma>0$ with $\gamma\in[0,1)$ \cite{Secchi2010}.

The basis for the Lyapunov--Schmitd type perturbation argument is the nondegeneracy of the groundstate of the limit equation \eqref{e-V-0},
which at present is not known beyond the slightly superquadratic case \(N = 3\), \(\alpha = 2\) and \(p > 2\) close to \(2\) \cite{Xiang},
as already discussed in \S\ref{s-nodegeneracy}.
An alternative approach to construct semiclassical solutions which does not rely on nondegeneracy are \emph{variational penalisation method}.
A penalisation of the \(L^p\) norm \cite{ByeonWang2003} was used to construct semiclassical multibump solutions to Choquard equations concentrating around minima \(V\), possibly in the presence of a magnetic field \cite{CingolaniSecchiSquassina2010}.

The existence of a global groundstate has been proved in the semiclassical limit when \(p \ge 2\),
\(V (x) = o (\abs{x}^\gamma)\) for some \(\gamma > 0\) and \(V\) stays away from \(0\) suitably
\cite{YangDing2013}. This results has been extended to similar problems in the presence of a magnetic field \citelist{\cite{YangWei2013}\cite{SunZhang2014}} and to corresponding \(p\)--Laplacian problems \citelist{\cite{YangWei2013}\cite{SunZhang2014}}.

Except of \cite{Secchi2010}, all of these results were obtained under the assumption $\liminf_{|x|\to\infty} V(x)>0$.
A penalisation method which allows to handle potentials without any restrictions on the decay (or growth) at infinity in the spirit of M.\thinspace{}del Pino and P.\thinspace{}Felmer \cite{delPino-Felmer-97} (see also \citelist{\cite{MorozVanSchaftingen2010}\cite{BonheureVanSchaftingen2008}}) was introduced in \cite{MorozVanSchaftingen2015CVAR}.

To clarify the issues related to the decay of the potential $V$, we consider the case when $V\ge 0$ and the support of \(V\) is a compact subset of $\Rset^N$ with $N\ge 3$.
This is of course the worst possible scenario from the point of view of the decay of $V$.
The rescaled semiclassical solutions $v_\eps$ of \eqref{e-V-eps-rescaled} then must satisfy for some $\rho>0$ the inequality at infinity,
\begin{equation}\label{e-V-eps-rescaled-infiniy}
-\Delta v_\eps\ge \bigl(I_\alpha\ast v_\eps^p\bigr)v_\eps^{p-1}\qquad\text{in $\Rset^N\setminus B_{\rho/\eps}$}.
\end{equation}
At the same time we expect $v_\eps$ to converge to the groundstates of the limit equation \eqref{e-V-0}.
While groundstates of \eqref{e-V-0} exist for $p\in\big(\frac{N+\alpha}{N},\frac{N+\alpha}{N-2}\big)$,
the existence of solutions to the inequality~\ref{e-V-eps-rescaled-infiniy} is restricted by the Liouville theorem~\ref{Thm-free}.
The admissible range for the existence of semi-classical solutions of \eqref{e-V-eps} is then given by the intersection of both regimes.

The \emph{penalised nonlinearity} \(g_\eps:\Rset^N\times\Rset\to \Rset\) defined in \cite{MorozVanSchaftingen2015CVAR} as
\[
 g_\varepsilon (x, s): = \chi_\Lambda(x) s_+^{p-1}+\chi_{\Rset^N\setminus\Lambda}(x)\min\big( s_+^{p-1},  H_\varepsilon (x)\big).
\]
where $\Lambda\subset\Rset^N$ is a potential well around the local minimum point $a$ of the potential $V$,
and $H_\eps$ is a \emph{penalisation potential}, which penalises large values of $s$ outside $\Lambda$.
In case $\liminf_{|x|\to\infty} V(x)>0$ the penalisation $H_\eps$ is usually chosen as a positive constant.
In \cite{MorozVanSchaftingen2015CVAR} the penalisation $H_\eps$ is constructed as
$$H_\eps(x)\approx c_\eps U(x)^{p-1},$$
where $c_\eps>0$ is a suitable constant and $U$ is a positive solution of inequality \eqref{e-V-eps-rescaled-infiniy}, provided that such a solution exists.
With such a choice of $H_\eps$ it is possible to show using the Stein--Weiss convolution inequality \eqref{eqHLSSW},
that the penalised functional
\[
\mathcal J_\eps(u)= \frac{1}{2}\int_{\Rset^N} \bigl( \eps^2\abs{\nabla u}^2 + V(x)\abs{u}^2\bigr)
 - \frac{p}{2 \eps^\alpha} \int_{\Rset^N} \bigabs{I_\frac{\alpha}{2}\ast G_\varepsilon (u)}^2,
\]
where $G_\eps(x,s)=\int_0^s g_\eps(x,t)\dif t$,
is well-defined and satisfies all the assumptions of the mountain-pass lemma in the space $H^1_V(\Rset^N)$.
Note that the action functional \eqref{eqAutoHomoAction} of the limit problem \eqref{e-V-0} is defined in $H^1(\Rset)\neq H^1_V(\Rset^N)$,
so the limit and perturbed problems are posed in the different function spaces!
%

A concentration analysis combined with the use of a nonlocal comparison principle allow to conclude in the case $p\ge 2$ that for $\eps > 0$ close enough to \(0\),
the mountain-pass critical point of the functional $\mathcal J_\eps$ solves the original problem \eqref{e-V-eps} and concentrates to a local minimum inside the set $\Lambda$.
The case $p<2$ remains open and we are not aware of any results about the existence of semiclassical solutions for \eqref{e-V-eps},
although $p<2$ includes some admissible regimes.
For $p\ge 2$ the results in \cite{MorozVanSchaftingen2015CVAR} are optimal from the point of view of the assumption on the admissible range of $p$ and decay of $V$.

\begin{theorem}[\cite{MorozVanSchaftingen2015CVAR}]
\label{theoremLinear}
Let \(\alpha \in ((N-4)_+, N)\) and \(V \in C(\Rset^N;[0, \infty))\).

$i)$ Assume that $p=2$ and either \(\alpha < N - 2\) or
\begin{equation}\label{e-positivity-quadratic}
 \inf_{x \in \Rset^N} V (x) \bigl(1 + \abs{x}^{N - \alpha}\bigr) > 0.
\end{equation}

$ii)$ Assume that \(p \in (2, \frac{N + \alpha}{(N - 2)_+})\) and either
$p > 1 + \frac{\max(\alpha, 1 + \frac{\alpha}{2})}{(N - 2)_+}$
or
\[
  \liminf_{\abs{x} \to \infty} V (x) \abs{x}^{2} > 0.
\]
If \(\Lambda \subset \Rset^N\) is an open bounded set such that
\[
  0 < \inf_{\Lambda} V < \inf_{\partial \Lambda} V,
\]
then the equation \eqref{e-V-eps} has a family of positive solutions \((u_\varepsilon)\) such that for
a family of points \((a_\varepsilon)\subset\Lambda\) and for every \(\rho > 0\), $\lim_{\varepsilon \to 0} V (a_\varepsilon)  = \inf_{\Lambda} V$ and
\begin{align*}
  \liminf_{\varepsilon \to 0} \varepsilon^{-N} \int_{B_{\varepsilon \rho} (a_\varepsilon)} \abs{u_\varepsilon}^2 &> 0,&
  \lim_{\substack{R \to \infty\\ \varepsilon \to 0}} \norm{u_\varepsilon}_{L^{\infty} (\Rset^N \setminus B_{\varepsilon R} (a_\varepsilon))}  &= 0.
\end{align*}
\end{theorem}

The global positivity assumption in \eqref{e-positivity-quadratic} cannot be removed.
In fact, an interesting essentially nonlocal phenomenon occurs in the range $\alpha> N-2$ for $p=2$.
If the potential $V(a_0)=o(|x-a_0|)^{\frac{4}{\alpha+2-N}-2}$ at some $a_0\in\Rset^N$ then for any solution $u_\eps$ of \eqref{e-V-eps}
and for any compact set $K\subset \Rset^N \setminus \{a_0\}$,
$$\int_Ku_\eps^2=o(\eps^N)\qquad\text{as \(\eps \to 0\).}$$
In particular, equation \eqref{e-V-eps} can not have solutions that concentrates inside $\Lambda$.
The nonlocal interaction with a zero of $V$ forces the rescaled solution to vanish everywhere outside this zero.
An even more delicate behaviour occurs in the case $\alpha=N-2$, see \cite{MorozVanSchaftingen2015CVAR}*{theorems 4 and 5}.

In the planar case \(N = 2\), semiclassical solutions have been constructed for an exponential nonhomogeneous nonlinearity \cite{AlvesCassaniTarsiYang}.
For \(N \ge 3\), similar results have been obtained for a general homogeneous nonlinearity as in \S\ref{sectionGeneralNonlinearity} \cite{YangZhangZhang2016}, following a strategy of \citelist{\cite{ByeonJeanjean2007a}\cite{ByeonJeanjean2007b}} for local problems.
Semiclassical solutions for Choquard equations involving $p$--Lpalacian are studied in \citelist{\cite{AlvesYang2014}\cite{AlvesYang2014JMP}}

\begin{problem}
Construct a concentrating family of solutions of \eqref{e-V-eps} in the subquadratic case \(p \in (\frac{N + \alpha}{N},2)\).
\end{problem}

Some existence results in this direction have been obtained when the external potentials $V$ is confining \cite{YangDing2013}.

\subsubsection{Strong electric field r\'egime}

When \(\frac{N - 2}{N + \alpha} < \frac{1}{p} < \frac{N}{N + \alpha}\), the problem
\[
 -\Delta u + (1 + \mu V_0)u = \bigl(I_\alpha \ast \abs{u}^p\bigr)\abs{u}^{p - 2}u,
\]
admits solutions for \(\mu > 0\) sufficiently large
if the function \(V_0\) is nonnegative, vanishes on an open set and has a positive sublevel set of finite measure; the solutions converge as \(\mu \to \infty\) to solutions of a Choquard problem on the zero set of the function \(V_0\) \cite{Lu2015}.
When the zero set of the function \(V_0\) has \(k\) connected components, then the problem has at least \(2^k - 1\) distinct nontrivial solutions when \(\mu\) is large enough \cite{AlvesNobregaYang}.

\subsubsection{$L^2$--constrained groundstates and symmetry breaking}

In \cite{AschbackerFrohlichGrafSchneeTroyer2002} the authors have considered the constrained minimisation problem
$$M_\lambda:=\inf\Big\{\mathcal{A}_{V}(u)\,:\,u\in H^1(\Rset^N),\,\int_{\Rset^N}u^2=\lambda\Big\},$$
in the case where $V\in C(\Rset^N)$ is a bounded external potential such that the bottom of the spectrum of $-\Delta+V$
is an eigenvalue.
It is proved that if $\lambda>0$ is sufficiently small then the minimisation problem $M_\lambda$ admits a unique positive minimiser,
while for large $\lambda$ minimisers of $M_\lambda$ concentrate around local minima of $V$.
In the case when there are several equivalent up to the symmetry local minima of $V$ this means that
a minimiser has to concentrate to one of them. As a consequence, the uniqueness and underlying symmetry of the minimisers
(wit respect to the symmetries of $V$) can not be preserved.

Related issues in the case of a confining potential $V$ have been discussed in \citelist{\cite{DengLuShuai2015}\cite{Ye2015arxiv}}.

\subsection{Magnetic potential}

For a given magnetic potential \(A : \Rset^N \to \Rset^N\) and an electric potential \(V : \Rset^N \to \Rset\), the Choquard equation with a variable electromagnetic field
imposes the function \(u : \Rset^N \to \Cset\) to satisfy
\begin{equation}
 (-i\nabla + A)^2 u + V u = \bigl(I_\alpha \ast \abs{u}^p\bigr) \abs{u}^{p - 2} u
 \qquad \text{in \(\Rset^N\)}.
\end{equation}
The \emph{constant magnetic fields and vanishing electric field case} was reviewed in \S\ref{sectionConstantMagneticField} above.

The \emph{existence of groundstates} has been proved when \(N = 3\), \(p = 2\) and \(\alpha = 2\) under trapping conditions, that is, when \(V\) is a reasonably bounded and localised perturbation of a constant and the groundstate level satisfies a strict inequality with a problem at infinity \cite{GriesemerHantschWellig2012}*{theorem~2.4}.

When the potentials \(V\) and \(A\) are compatible with the action of a group of isometries on \(\Rset^N\), \emph{infinitely many solutions} have been constructed when all the orbits of points in \(\Rset^N \setminus \{0\}\) are infinite \cite{CingolaniClappSecchi2012}*{theorem~1.1} and the existence of several solutions has been proved under weaker assumption on the orbits but more stringent assumptions on the magnetic potential \(A\), implying its decay at infinity
\citelist{%
\cite{CingolaniClappSecchi2012}*{theorem~1.2}%
\cite{Salazar2015}%
}

The study of the magnetic Choquard equation had been initiated earlier by S. Cingolani, S. Secchi and M. Squassina \cite{CingolaniSecchiSquassina2010}, with the analysis of the \emph{semi-classical problem}
\[
 (-i\varepsilon \nabla + A)^2 u + V u = \bigl(I_\alpha \ast \abs{u}^p\bigr) \abs{u}^{p - 2} u
 \qquad \text{in \(\Rset^N\)}.
\]
For \(V\) positive, solutions for \(\eps > 0\) small enough concentrating around finitely many local minimum points as \(\eps \to 0\) have been constructed \cite{CingolaniClappSecchi2013}*{theorem~1.1}.
Under symmetry assumptions on the electric potential \(V\) and the magnetic potential \(A\),
multiple families of solutions concentrating around some points have been constructed
\citelist{
\cite{YangWei2013}
\cite{Salazar2015}
}.
Multiple concentrating families of solutions exist when \(\alpha > N - 2\) with a nonhomogeneous nonlinearity \cite{AlvesFigueiredoYang2016}.
When \(V\) achieves \(0\) as a minimum at one point of \(\Rset^N\), solutions concentrating around that point have been constructed \cite{YangWei2013}.

\begin{problem}
Understand the semi-classical limit in the strong magnetic field r\'egime \(A \simeq \varepsilon^{-1}\), as it has been performed for the nonlinear Schr\"odinger equation \cite{DiCosmoVanSchaftingen2015}.
\end{problem}

\subsection{Nonlinear perturbation and other equations}


\subsubsection{Source term}
The Choquard equation
\begin{equation}
 \label{eqChoquard-Vf}
 -\Delta u + Vu = \bigl(I_\alpha \ast \abs{u}^p\bigr)|u|^{p-2}u+\mu f(x) \quad \text{in \(\Rset^3\)},
\end{equation}
with $N=3$, $\alpha=2$, $p=2$, $V\equiv 1$ and a source term $f\in H^{-1}(\Rset^3)$ has been studied in \cite{KupperZhangXia2003}: for some $\mu^{**}\ge \mu^*>0$ equation \eqref{eqChoquard-Vf}
possesses at least two positive solutions for $\mu\in(0,\mu^*)$ and no positive solution for $\mu>\mu^{**}$,
in accordance with an earlier result for local nonlinear equations on bounded domains \cite{Tarantello1992} (see also \citelist{\cite{Zhang1996}\cite{Zhang2000}\cite{Zhang2001}}).
A similar result in the case of a confining potential $V$ for $N\ge 3$, \(\alpha \in (0, N)\) and \(\frac{N - 2}{N + \alpha} < \frac{1}{p} < \frac{N}{N + \alpha}\) has been obtained \cite{XieXiaoWang2015}.

\subsubsection{Nonautonomous nonlocal term}
The existence of groundstates for the nonautonomous equation
\[
 -\Delta u + u = \bigl(I_\alpha \ast K \abs{u}^{p}\bigr)K \abs{u}^{p - 2} u,
\]
where \(K \in C (\Rset^N)\) is a nonnegative nonlinear potential and \(\frac{N - 2}{N + \alpha} < \frac{1}{p} < \frac{N}{N + \alpha}\)
has been studied \cite{Wang2016}.
Other nonautonomous modifications of the nonlocal term have been considered
\citelist{%
\cite{ZhangKupperHuXia2006}%
\cite{Cao1989}%
}.

\subsubsection{Local nonautonomous perturbation}
A number of authors considered nonautonomous perturbation of the Choquard equation of the type
\begin{equation}\label{eChoqLocalq}
 -\Delta u + V u = \bigl(I_2 \ast K\abs{u}^p\bigr)K|u|^{p-2}u+a|u|^{q-2}u
 \qquad \text{in \(\Rset^3\)}.
\end{equation}
The corresponding autonomous case was discussed in \S\ref{sectLocalQ}.

In \cite{Vaira2011}, the case $p=2$, $q\in(2,6)$ and $V\equiv 1$ has been considered.
The potentials $K,a:\Rset^3\to\Rset$ are positive functions such that $K(x)-K_\infty\in L^2(\Rset^3)$ and $a(x)-a_\infty\in L^\frac{6}{6-q}(\Rset^3)$, for some positive constants $K_\infty$, $a_\infty$.
Under additional structural conditions on $K$ and $a$, \cite{Vaira2011} has established the existence of a positive groundstate for \eqref{eChoqLocalq}.
In \cite{Vaira2013}, under additional assumption $q=3$ the existence of a positive bound state has been proved in the situation when groundstates may not exists. The paper also addresses the uniqueness and nondegeneracy of the radial groundstate
of the associated autonomous limit problem
\[
 -\Delta u + u = \bigl(I_2 \ast K_\infty\abs{u}^2\bigr)K_\infty u+a_\infty|u|u
 \qquad \text{in \(\Rset^3\)},
\]
see \cite{Vaira2013}*{proposition~3.6 and theorem 3.7}.

Nonautonomous Choquard equations of type \eqref{eChoqLocalq} with $p=2$, constant potentials $V$ and $K$ and a general local nonlinearity $f(x,u)$
superlinear subcritical type instead of $a|u|^{q-2}u$ have been considered in \citelist{\cite{Mugnai2011}\cite{JeongSeok2014}}
(for related systems of equations, see \citelist{\cite{Hajaiej2012}\cite{Hajaiej2013}}).
Equation \eqref{eChoqLocalq} with a general $p$ and $q$ on bounded domains has been studied in \cite{AzzolinidAveniaLuisi2013}.

Equation \eqref{eChoqLocalq} with $p=2$, $q\in(3,6)$, constant potentials $K>0$ and $a<0$, and a periodic potential $V$ such that $0$ is in a spectral gap of $-\Delta+V$ has been studied in \cite{ChenXiao2015}.
The case $q=\frac{10}{3}$ is related to a Hartree type model for crystals \citelist{\cite{CattoLeBrisLions2002}\cite{LeBrisLions2005}}.
The existence of a nontrivial solution has been proved in \cite{ChenXiao2015} for all sufficiently small $K$ using a linking theorem.

Groundstates of the equation \eqref{eChoqLocalq} when \(K\) is a perturbation of \(\abs{x}^{-\beta}\) with \(\beta \in (0, 2 - \frac{N - \alpha}{2})\) and \(p = 2\) in \(\Rset^N\) with \(N \ge 3\)
have been constructed \cite{GuoSu2016}.


\appendix
\section{About the Riesz potentials}
\label{sectionAppendix}

The Riesz potential have been introduced by M.\thinspace{}Riesz in the 1930s \cite{Riesz1939} and was systematically studied in his fundamental paper \cite{Riesz1949}.
An exposition on basic functional--analytic properties of the Riesz potentials could be found in E.\thinspace{}Stein's monograph \cite{Stein1970}*{\S 5.1} and also in many places in \cite{LiebLoss2001}.
A systematic potential theoretic study of the Riesz potentials is presented in the monographs by N.\thinspace{}Landkof \cite{Landkof}
(see also N.\thinspace{}du\thinspace{}Plessis \cite{duPlessis}*{chapter 3} for a shorter exposition).

\subsection{Definition and semigroup property}
The Riesz potential \(I_\alpha\) of order \(\alpha \in (0, N)\) on the Euclidean space \(\Rset^N\) of dimension \(N \ge 1\) is defined for each \(x \in \Rset^N \setminus \{0\}\) by
\begin{equation*}
I_\alpha (x) = \frac{A_\alpha}{\abs{x}^{N - \alpha}},
\qquad\text{where $A_\alpha = \frac{\Gamma(\tfrac{N-\alpha}{2})}{\Gamma(\tfrac{\alpha}{2})\pi^{N/2}2^{\alpha}}$.}
\end{equation*}
The choice of normalisation constant $A_\alpha$ ensures that
the \emph{semigroup property}
\[
  \label{eqRieszSemigroup}
  I_\alpha\ast I_\beta=I_{\alpha+\beta},\qquad\text{$\forall\alpha,\beta>0$ such that $\alpha+\beta<N$,}
\]
and, for $N\ge 3$, the property
$$-\Delta I_\alpha=I_{\alpha-2},\qquad\forall\alpha\in(2,N).$$
In addition, $-\Delta I_2=\delta$, where $\delta$ is the Dirac delta function, that is $I_2$ is the Green function of the Laplacian $-\Delta$ on $\Rset^N$.
More generally, $I_\alpha$ could be interpreted as the inverse of the fractional Laplacian operator $(-\Delta)^{\alpha/2}$.
See e.g. \citelist{\cite{Riesz1949}\cite{Landkof}*{\S I.1}\cite{Stein1970}*{\S V.1.1}} for the study of these fundamental properties of the Riesz potentials.

When $\alpha\to 0$, $I_\alpha\to\delta$, in the vague sense \cite{Landkof}*{p.\thinspace{}46}.
When $\alpha\to N$, $I_\alpha*f\to A_N\log\big(\frac{1}{|x|}\big)*f$ for every $f\in C^\infty_c(\Rset^N)$ such that $\int_{\Rset^N}f=0$, where $A_N=\lim_{\alpha\to N}(N-\alpha)A_\alpha=1/\big(\Gamma(d/2)\pi^{d/2}2^{d-1}\big)$ \cite{Landkof}*{p.\thinspace{}50}.

The definition of the Riesz potentials as well as the semigroup property could be extended from $\alpha\in(0,N)$ to arbitrary complex $\alpha$
with $\mathrm{Re}(\alpha)>0$ and $\frac{\alpha-N}{2}\not\in\Nset\cup\{0\}$, however then the convolution $I_\alpha\ast f$ should be interpreted in the distributional sense, see for example \cite{Landkof} or \cite{Samko02}*{chapter 2} for a more recent exposition.
In this survey we always assume that $\alpha\in(0,N)$ and the convolution $I_\alpha\ast f$ is understood in the sense of the Lebesgue integral.

\subsection{$L^p$--estimates}
The Riesz potential of order $\alpha\in(0,N)$ of a function $f\in L^1_{\mathrm{loc}}(\Rset^N)$ is defined as
$$I_\alpha\ast f(x):=A_\alpha\int_{\Rset^N}\frac{f(y)}{\abs{x-y}^{N - \alpha}}dy.$$
The latter integral converges in the classical Lebesgue sense for a.e. $x\in\Rset^N$ if and only if
\begin{equation}\label{e-Riesz-necessary}
f\in L^1\big(\Rset^N,(1+|x|)^{-(N-\alpha)}dx\big),
\end{equation}
Moreover, if \eqref{e-Riesz-necessary} does not hold then $I_\alpha\ast|f|=+\infty$ everywhere in \(\Rset^N\) \cite{Landkof}*{p.61-62}.

The Riesz potential $I_\alpha$ is well--defined as an operator on the whole space $L^q(\Rset^N)$ if and only if $q\in[1,\frac{N}{\alpha})$.
The Hardy--Littlewood--Sobolev inequality \citelist{\cite{HardyLittlewoodPolya1952}*{theorem~382}\cite{Sobolev1938}}
(see also \citelist{\cite{Stein1970}*{theorem~V.1}\cite{LiebLoss2001}*{theorem~4.3}}),
which states that if \(q \in (1, \infty)\)
and if \(\alpha < \frac{N}{q}\), then for every \(f \in L^q (\Rset^N)\), we have \(I_\alpha \ast f \in L^{\frac{Nq}{N - \alpha q}} (\Rset^N)\)
and
\begin{equation}
\label{eqHLS}
 \Bigl(\int_{\Rset^N} \bigabs{I_\alpha \ast f}^\frac{Nq}{N - \alpha q}\Bigr)^{\frac{1}{q} - \frac{\alpha}{N}}
 \le C_{N, \alpha, q} \Bigl(\int_{\Rset^N} \abs{f}^q\Bigr)^\frac{1}{q}.
\end{equation}
If \(q=\frac{2N}{N + \alpha}\), then the optimal constant is given by \citelist{\cite{Lieb1983}\cite{FrankLieb2010}\cite{LiebLoss2001}*{theorem~4.3}}
\begin{equation}
  C_{N, \alpha, q} = \frac{\Gamma (\frac{N - \alpha}{2})}{2^\alpha \pi^{\alpha/2} \Gamma (\frac{N + \alpha}{2}) }
  \biggl(\frac{\Gamma(\frac{N}{2})}{\Gamma (N)}\biggr)^\frac{\alpha}{N}.
\end{equation}
The inequality \eqref{eqHLS} implies that if \(q \in (1, \infty)\), \(\alpha < \frac{N}{q}\) and $\frac{1}{r}=\frac{1}{q}-\frac{\alpha}{N}$, then
$$I_\alpha: L^q (\Rset^N)\to L^{r} (\Rset^N)$$
is a bounded linear operator.
If $f\in L^1(\Rset^N)$ then in general $I_\alpha*f\not\in L^\frac{N}{N-\alpha}(\Rset^N)$ (see for example \cite{Stein1970}*{\S V.1.2}),
however
\[
  I_\alpha: L^1 (\Rset^N)\to L^{r}\big(\Rset^N,(1+|x|)^{-\lambda}dx\big)
\]
is a bounded operator for any $r\in[1,\frac{N}{N-\alpha})$ and $\lambda>N-r(N-\alpha)$ \cite{Samko02}*{p.\thinspace{}38}.
When $q\ge\frac{N}{\alpha}$ then $I_\alpha$ is not well-defined on the whole space $L^q(\Rset^N)$.
However, if $f\in L^{\frac{N}{\alpha}}(\Rset^N)$ and we additionally assume that $I_\alpha*f$ is almost everywhere finite on $\Rset^N$ then $I_\alpha*f$
is a function of bounded mean oscillation (BMO) \citelist{\cite{SteinZygmund1967}*{theorem~2}\cite{MuckenhouptWheeden1974}*{theorem~7}}.
If $I_\alpha*f$ is almost everywhere finite on $\Rset^N$, $f\in L^{q}(\Rset^N)$ and $\frac{N}{q}<\alpha<\frac{N}{q}+1$ then
$I_\alpha*f$ is H\"older continuous of order $\alpha-\frac{N}{q}$ \cite{duPlessis1955}*{theorem~2}.

These mapping properties of the Riesz potentials are important not only as a tool to control the domain of definition
of the action functional of the Choquard equation but also in the study of the regularity properties of solutions by bootstrap type procedures
(see \S\ref{sectRegularity}).

The weighted version of the Hardy--Littlewood--Sobolev inequality due to E. Stein and G. Weiss \cite{SteinWeiss1958} states that
if \(q \in (1, \infty)\), $s<N (1 - \frac{1}{q})$, $t<\frac{N}{r}$, $s+t\ge 0$, $q\le r<\infty$ and $\frac{1}{r}=\frac{1}{q}+\frac{s+t-\alpha}{N}$
then for any $f\in L^q(\Rset^N,|x|^{sq}dx)$,
\begin{equation}
\label{eqHLSSW}
 \Bigl(\int_{\Rset^N} \bigabs{I_\alpha \ast f(x)}^r|x|^{-rt} dx\Bigr)^{\frac{1}{r}}
 \le C \Bigl(\int_{\Rset^N} \abs{f(x)}^q|x|^{sq} dx\Bigr)^\frac{1}{q}.
\end{equation}
When $r=q$ the inequality is also known as the Hardy--Rellich inequality.
Optimal constants for \eqref{eqHLSSW} are known in some special cases, see \citelist{\cite{He1977}\cite{Yafaev1999}\cite{Samko05}\cite{Beckner2008}
  \cite{Beckner2012}
  \cite{MorozVanSchaftingen2012}
  }.

B. Rubin \cite{Rubin-1982} has proved that for radial functions $f\in L^q_{\mathrm{rad}}(\Rset^N,|x|^{sq}dx)$ the same inequality \eqref{eqHLSSW} holds for
a wider range $s+t\ge -(N-1)\big(\frac{1}{q}-\frac{1}{r}\big)$ (see also \citelist{\cite{Duoandikoetxea13}\cite{DeNapoli-2014}}).
Weighted inequalities of Stein--Weiss type become important in the analysis of nonautonomous Choquard equations with confining
(\S\ref{sectConfining}) or decaying (\S\ref{sectSemiclass}) potentials.

\subsection{Energy properties}
The Riesz potential \(I_\alpha\) naturally induces the quadratic form \(D_\alpha\) defined by
\[
  D_{\alpha}(f,g):=A_\alpha\iint_{\Rset^N \times \Rset^N}\frac{f(x)g(y)}{\hphantom{!}\abs{x-y}^{N - \alpha}}\dif y\dif x.
\]
A direct consequence of the semigroup property \eqref{eqRieszSemigroup} is the inequality
\begin{equation}
  D_{\alpha}(f,f)
  =\int_{\Rset^N} \bigl(I_\alpha \ast f) f
  =\int_{\Rset^N}\big|I_{\alpha/2}*f|^2\ge 0,
\end{equation}
valid for all functions $f$ such that $D_\alpha(|f|,|f|)<\infty$.
Moreover, $D_\alpha(f,f)=0$ if and only if $f\equiv 0$ \cite{Landkof}*{theorem~1.15}.

The Riesz--Sobolev rearrangement inequality states that for any two nonnegative functions $f,g$ such that $D(f,g)<\infty$,
\begin{equation}%
\label{eqRearrangement}%
D_\alpha(f^*,g^*)\le D_\alpha(f,f),
\end{equation}
where $f^*$ denotes the symmetric decreasing rearrangement of $f$.
It was first established by F.\thinspace{} Riesz \cite{RieszF} in one dimension, then S.\thinspace{}L.\thinspace{}Sobolev \cite{Sobolev1938} extended the result to $\Rset^N$ (see also \cite{BrascampLiebLuttinger1974}).
The equality in \eqref{eqRearrangement} occurs if and only if \(f\) is the translation of radially symmetric and nonincreasing function \citelist{\cite{Lieb1977}*{lemma 3}\cite{Burchard1996}\cite{LiebLoss2001}*{\S 3.7--3.9}}.
Inequality \eqref{eqRearrangement} is fundamental in the study of radial symmetry of groundstates of Choquard equation (\S\ref{sectSymmetry}).

\subsection{Positivity and decay estimates}
When the function $f\in L^1_{\mathrm{loc}}(\Rset^N)$ is nonnegative, an elementary estimate shows that for every $x\in\Rset^N$,
$$I_\alpha*f(x)\ge\frac{A_\alpha}{R^{N-\alpha}}\int_{B_R(x)}f(y)\dif y.$$
In particular, if the function \(f\) is positive on a set of positive measure of \(\Rset^N\), then $I_\alpha*f$ is everywhere strictly positive on $\Rset^N$.
Similarly, for each $x\in\Rset^N$ one can estimate
\begin{equation}\label{eqLowerbound}
I_\alpha*f(x)\ge\frac{A_\alpha}{(2|x|)^{N-\alpha}}\int_{B_{2|x|}(x)}f(y)\dif y\ge \frac{c}{|x|^{N-\alpha}},
\end{equation}
that is $I_\alpha*f$ can not decay faster then $I_\alpha$ at infinity, even if the function $f$ is compactly supported.

These decay properties of the Riesz potential are essential for the study of asymptotic decay of groundstates
of Choquard equations (\S\ref{s-Decay}) and Liouville's theorems (\S\ref{sLiouville}).
To illustrate the decay of the Riesz potentials at infinity, assume that the pointwise bound
\begin{equation*}\label{eqgamma}
\limsup_{\abs{x} \to \infty} f(x) \abs{x}^\gamma <\infty
\end{equation*}
holds.
Then, by a direct computation \cite{MorozVanSchaftingen2013JDE}*{lemma A.1},
\begin{align*}
\limsup_{\abs{x} \to \infty}\, (I_\alpha \ast f)(x)\,\abs{x}^{\gamma-\alpha} &<\infty& & \text{if \(\alpha<\gamma<N\)},\\
\limsup_{\abs{x} \to \infty}\, (I_\alpha \ast f)(x) \frac{\hphantom{l}\abs{x}^{N - \alpha}}{\log \abs{x}} &<\infty&
&\text{if \(\gamma=N\)},\bigskip\\
\limsup_{\abs{x} \to \infty}\, (I_\alpha \ast f)(x)\,\abs{x}^{N - \alpha} &<\infty& & \text{if \(\gamma>N\)},
\end{align*}
and the bounds are optimal, as can be seen by choosing $f(x)=(1+|x|)^{-\gamma}$ in the case $\gamma\ge N$,
and by comparing with \eqref{eqLowerbound} in the case $\gamma>N$.

If $\gamma>N$ then the decay of $I_\alpha*f$ explicitly depends on $\|f\|_{L^1(\Rset^N)}$.
More specifically \cite{MorozVanSchaftingen2013JFA}*{lemma 6.2}, assume that
\begin{equation}\label{eqgammaN}
\sup_{\Rset^N}|f(x)|(1+|x|)^\gamma<\infty,
\end{equation}
for some $\gamma>N$. Then
\begin{equation}\label{eqRieszLocalised}
I_\alpha*f(x)=\Big(I_\alpha(x)\int_{\Rset^N} f (y) \dif y\Big)\big(1+o(1)\big)\qquad\text{as $|x|\to\infty$}.
\end{equation}
Note that the assumption $f\in L^1(\Rset^N)$ alone does not imply that $I_\alpha*f=O(|x|^{-(N-\alpha)})$ even if $f$ is radial:
some additional control on the decay of $f$ at infinity is always needed \cite{SiegelTalvila99}.
However, if $f$ is radial and $\alpha>1$ then $I_\alpha*f=O(|x|^{-(N-\alpha)})$ if and only if $f\in L^1_{\mathrm{rad}}(\Rset^N)$ \cite{SiegelTalvila99}*{theorem~5.A}.
Radial estimates on the Riesz potentials could be found in \citelist{\cite{Rubin-1982}\cite{Duoandikoetxea13}\cite{DeNapoli-2014}\cite{Thim2016}\cite{MercuriMorozVanSchaftingen}}.


\end{document}